\documentclass[12pt,reqno]{amsart}
\usepackage{amsmath,amssymb,amsthm}
\usepackage[shortlabels]{enumitem}
\usepackage{tikz}
\usepackage{hyperref}
\usepackage[sort]{cite}
\usepackage{booktabs,tabularx}
 \usepackage{booktabs,tabularx}
 
\newtheorem{theorem}{Theorem}[section]
\newtheorem{lemma}[theorem]{Lemma}
\newtheorem{proposition}[theorem]{Proposition}
\newtheorem{corollary}[theorem]{Corollary}
\newtheorem{conjecture}[theorem]{Conjecture}
\newtheorem{question}[theorem]{Question}
\theoremstyle{definition}
\newtheorem{definition}[theorem]{Definition}
\newtheorem{example}[theorem]{Example}
\theoremstyle{remark}
\newtheorem{remark}[theorem]{Remark}

\newcommand{\kk}{\Bbbk}
\newcommand{\T}{\mathcal{T}}
\newcommand{\Ht}{\widetilde{H}}
\newcommand{\Sk}{\mathrm{B}}
\newcommand{\sym}[2]{#1^{\{#2\}}}
\newcommand{\hKL}{h_{\mathrm{KL}}}
\newcommand{\cH}{\mathcal{H}}
\newcommand{\cK}{\mathcal{K}}
\newcommand{\cS}{\mathcal{S}}
\newcommand{\cL}{\mathcal{L}}

\newcommand{\fF}{\mathfrak{F}}
\newcommand{\sqf}{\mathrm{sqf}}
\DeclareMathOperator{\reg}{reg}
\DeclareMathOperator{\pd}{pd}
\DeclareMathOperator{\Tor}{Tor}
\DeclareMathOperator{\hgt}{ht}
\DeclareMathOperator{\lk}{lk}
\DeclareMathOperator{\Ind}{Ind}
\DeclareMathOperator{\Cl}{Cl}
\DeclareMathOperator{\topp}{top}
\DeclareMathOperator{\cost}{cost}

\title[Tolerance complexes and squarefree symbolic powers]
      {Tolerance complexes of chordal graphs and squarefree symbolic powers}
\author{S Selvaraja}
\email{selvas@iitbbs.ac.in}
\address{Department of Mathematics,
Indian Institute of Technology Bhubaneswar,
Bhubaneswar-752050, India}



\begin{document}

\begin{abstract}
For a simplicial complex $K$ and an integer $t\ge0$, the $t$-tolerance
complex $\mathcal{T}_t(K)$ consists of the sets that become faces of $K$
after deleting at most $t$ of their elements. Kim and Lew conjectured that
if $K$ is $d$-Leray, then $\mathcal{T}_t(K)$ is $(\eta(d+1,t+1)-1)$-Leray
(the Leray conjecture), and that if $K$ is $d$-collapsible, then
$\mathcal{T}_t(K)$ is $(\eta(d+1,t+1)-1)$-collapsible (the collapsibility
conjecture), where $\eta$ is the Erd\H{o}s--Gallai number.
We prove the collapsibility conjecture for $d=1$ when $t\le3$ or
$\dim K\le4$, and for all $t$ when $K$ is the clique complex of a subtree
graph of a spider, a class of chordal graphs containing all interval and
split graphs; in particular, it holds for nerves of finite families of
intervals. For every $d$, we prove both conjectures when $K$ has at most
$d+1$ facets, when $K$ has nonzero homology in its top dimension, and when
the Stanley--Reisner ideal of $K$ is squarefree stable. 
The Stanley--Reisner ideal of $\mathcal{T}_t(K)$ is the $(t+1)$-st
squarefree symbolic power of that of $K$, so the Leray conjecture becomes a
regularity bound. Using this, we prove the Leray conjecture for flag
complexes when $t\le3$, and for all $t$ for independence complexes of
bipartite and chordal graphs, and we show that it is preserved under joins.
We also construct chordal graphs $H$ for which the Leray number of
$\mathcal{T}_t(\mathrm{Cl}(H))$ exceeds its Helly dimension by an
arbitrarily large amount, and others for which it depends on the field.
Finally, we reduce the case $t=1$ of the Leray bound for $d$-collapsible
complexes to a question on intersections of links.
\end{abstract}

\maketitle

\section{Introduction}\label{sec:intro}

Helly's theorem states that a finite family of convex sets in $\mathbb{R}^d$
has nonempty intersection whenever every subfamily of at most $d+1$ members
does. Much of the theory of intersection patterns of convex sets can be
studied through the \emph{nerve} of a family, the simplicial complex
recording its intersecting subfamilies. Wegner \cite{Weg75} proved that
nerves of finite families of convex sets in $\mathbb{R}^d$ are
$d$-collapsible and that every $d$-collapsible complex is $d$-Leray. Here a
complex $K$ is \emph{$d$-collapsible} if it can be reduced to the void
complex by successively removing faces of cardinality at most $d$, each
contained in a unique maximal face, together with all faces containing them.
It is \emph{$d$-Leray} over a field $\kk$ if every induced subcomplex has
vanishing reduced homology over $\kk$ in all dimensions at least $d$. The
least such integers are the \emph{collapsibility number} $C(K)$ and the
\emph{Leray number} $L_\kk(K)$. The \emph{Helly dimension} $h(K)$ is the
maximum of zero and the dimensions of the minimal nonfaces of $K$. These
invariants satisfy
$h(K)\le L_\kk(K)\le C(K).$
Helly's theorem, together with its fractional and colorful versions, extends
to $d$-Leray complexes, so bounds on Leray numbers yield Helly-type
theorems; see the survey \cite{Tan13}.

This paper concerns tolerant versions of these results. A family of sets has
a nonempty intersection \emph{with tolerance $t$} if its members have a
common point after the removal of at most $t$ of them. For a simplicial
complex $K$ on a finite set $W$ and an integer $t\ge0$, the
\emph{$t$-tolerance complex} is
\[
\T_t(K)=\{\sigma\subseteq W:\ \sigma\setminus A\in K
\text{ for some }A\subseteq\sigma \text{ with } |A|\le t\}.
\]
If $K$ is the nerve of a family of sets, then the faces of $\T_t(K)$ are the
subfamilies that become intersecting after removing at most $t$ members.
Eckhoff \cite{Eck03} and, independently, Montejano and Oliveros \cite{MO11}
proved a tolerant Helly theorem. In combinatorial form, it states that if
$h(K)\le d$, then $h(\T_t(K))\le\eta(d+1,t+1)-1$; see
\cite[Theorem~1.2]{KL23}. Here $\eta(r,k)$ is the
\emph{Erd\H{o}s--Gallai number}: the maximum number of vertices of an
$r$-uniform hypergraph with no isolated vertices and transversal number $k$
in which deleting any edge decreases the transversal number. Thus the effect
of tolerance is governed by a classical extremal parameter of hypergraph
theory. Geometric assumptions can improve this bound: for families of convex
sets in $\mathbb{R}^2$ with $t=3$, it suffices to check subfamilies of size
$14$, whereas $\eta(3,4)=15$ \cite{JKL21b}; see \cite[Section~5]{JKL21}.

Since $h\le L_\kk\le C$, bounds on Leray and collapsibility numbers are stronger than bounds on Helly dimensions, it is
natural to study the Leray and collapsibility numbers of tolerance
complexes. Kim and Lew \cite{KL23} proved that if $K$ is $d$-collapsible,
then $\T_t(K)$ is $\hKL(t,d)$-Leray, where $\hKL$ is an explicit function
with $\hKL(t,1)=2t+1$, and that $\T_1(K)$ is $5$-Leray whenever $K$ is
$2$-collapsible; they used these bounds to obtain tolerant colorful Helly
theorems. The values $\hKL(t,1)=2t+1$ and $5$ coincide with the Erd\H{o}s--Gallai bounds $\eta(2,t+1)-1$ and $\eta(3,2)-1$, respectively,
which led them to the following conjectures
\cite[Conjectures~1.3 and~1.4]{KL23}.

\begin{conjecture}[Kim--Lew \cite{KL23}]\label{conj:leray}
Let $\kk$ be a field and let $K$ be a $d$-Leray complex over $\kk$. Then, for
every integer $t\ge0$, the complex $\T_t(K)$ is $(\eta(d+1,t+1)-1)$-Leray
over $\kk$.
\end{conjecture}

\begin{conjecture}[Kim--Lew \cite{KL23}]\label{conj:coll}
Let $K$ be a $d$-collapsible complex. Then, for every integer $t\ge0$, the
complex $\T_t(K)$ is $(\eta(d+1,t+1)-1)$-collapsible.
\end{conjecture}

The proposed bounds cannot be improved when $d=1$ or $t=1$; see
\cite{KL23} for $d=1$ and the construction of Montejano and Oliveros
\cite{MO11} for $t=1$. Both conjectures remain open in general. Even the
weaker statement that $\T_t(K)$ is $(\eta(d+1,t+1)-1)$-Leray whenever $K$ is
$d$-collapsible is open, already for $t=1$ and $d=3$; for $t=1$ the cases $d\le2$ were settled by Kim and Lew \cite[Theorems~1.5 and~1.6]{KL23}.

The Leray conjecture also has a natural algebraic interpretation. Let
$S=\kk[x_w:w\in W]$ and let $I_K$ be the Stanley--Reisner ideal of $K$.
We observe that the Stanley--Reisner ideal of $\T_t(K)$ is the
\emph{squarefree symbolic power} $\sym{I_K}{t+1}$, the ideal generated by the
squarefree monomials in the symbolic power $I_K^{(t+1)}$
(Proposition~\ref{prop:sym}); for matroid complexes this was observed in \cite{KL24,MN24}, see also \cite{DFKT24}. Consequently,
$L_\kk(\T_t(K))=\reg(S/\sym{I_K}{t+1})$, and
Conjecture~\ref{conj:leray} becomes a bound on the regularity of these
quotients in terms of $\reg(S/I_K)$. The regularity of symbolic and
squarefree powers has been studied extensively
\cite{EHHS22,MNPTV22,SF23,CDRS}, and this correspondence allows results to be
transferred in both directions.
We first consider the case $d=1$ of Conjecture~\ref{conj:coll}. Since
$\eta(2,t+1)=2t+2$, it asserts that $\T_t(K)$ is $(2t+1)$-collapsible
whenever $K$ is $1$-collapsible. A nonvoid complex is $1$-collapsible if and
only if, on its vertex set, it is the clique complex $\Cl(H)$ of a chordal
graph $H$ (Lemma~\ref{lem:chordal}). A \emph{spider} is a tree with at most
one vertex of degree at least $3$, and a \emph{subtree graph of a spider} is
the intersection graph of a family of subtrees of a spider. These graphs are
chordal and include all interval and split graphs.

\begin{theorem}[Theorems~\ref{thm:spider} and~\ref{thm:smallt}]
\label{main:chordal}
Let $H$ be a chordal graph and let $t\ge0$. If $t\le3$ or $H$ is a subtree
graph of a spider, then $\T_t(\Cl(H))$ is $(2t+1)$-collapsible.
\end{theorem}

In particular, the case $d=1$ of Conjecture~\ref{conj:coll} holds for
$t\le3$, and for all $t$ for nerves of finite families of intervals. The corresponding Leray bound $L_\kk(\T_t(K))\le2t+1$ for $1$-collapsible $K$ was proved by Kim and Lew \cite[Theorem~1.5]{KL23}; the new content of Theorem~\ref{main:chordal} is collapsibility. The
proof is a deletion--link induction on a larger class of complexes,
controlled by a potential function.
For $t=0$, and for threshold graphs for all $t$, the three invariants of
$\T_t(\Cl(H))$ coincide, so one might expect the regularity of the
corresponding squarefree symbolic powers to equal the largest degree of a
minimal generator. The next result shows that this fails, with arbitrarily
large gaps.

\begin{theorem}[Theorem~\ref{thm:notHelly}]\label{main:notHelly}
For all integers $r,c\ge2$, there is a chordal graph $H$ such that, for
$t=c^r-1$ and every field $\kk$,
\[
h(\T_t(\Cl(H)))=c^r+c-1
~\text{and}~
L_\kk(\T_t(\Cl(H)))=c^r+r(c-1).
\]
Moreover, for every $t\ge3$, there is a chordal graph $H$ such that
$L_\kk(\T_t(\Cl(H)))>h(\T_t(\Cl(H)))$.
\end{theorem}

This is a special case of a general construction
(Theorem~\ref{thm:splitgap}) that relates any complex with $t+1$ facets to
the $t$-tolerance complex of a split graph. The same construction shows that
$L_\kk(\T_t(\Cl(H)))$ can depend on the field, so $C$ and $L_\kk$ differ in
general for chordal graphs (Corollary~\ref{cor:field}), although the bound
$C(\T_t(\Cl(H)))\le2t+1$ may still hold. For $t=1$, all three invariants
coincide (Proposition~\ref{prop:t1}).

For flag complexes, comparing $\sym{I}{k}$ with the symbolic power $I^{(k)}$
(Proposition~\ref{prop:compare}) and using known regularity bounds for
symbolic powers of edge ideals, we obtain the following.

\begin{theorem}[Corollary~\ref{cor:flag}]\label{main:flag}
Let $K$ be a flag complex, and suppose that $t\le3$, or that $K$ is the
independence complex of a bipartite or chordal graph, or of a graph with no
odd cycle of length at most $2t-1$. Then Conjecture~\ref{conj:leray} holds
for $K$ and $t$. Moreover, $L_\kk(\T_t(K))\le L_\kk(K)+2t$.
\end{theorem}

The Leray conjecture is also preserved under joins. The proof combines a
known regularity bound for sums of ideals in disjoint sets of variables
(Theorem~\ref{thm:sum}) with a superadditivity property of Erd\H{o}s--Gallai
numbers (Lemma~\ref{lem:etamerge}).

\begin{theorem}[Corollary~\ref{cor:joins}]\label{main:joins}
If $K_1$ and $K_2$ are complexes on disjoint sets that satisfy
Conjecture~\ref{conj:leray} over $\kk$, then so does their join $K_1*K_2$.
\end{theorem}

Since $\eta(d+1,t+1)\ge(d+1)(t+1)$, the linear estimates
$L_\kk(\T_t(K))\le L_\kk(K)+t$ and $C(\T_t(K))\le C(K)+t$ imply the
conjectured bounds. These estimates fail in general
(Remark~\ref{rem:beyondshifted}), but they hold for the following classes,
often with equality. For a matroid $M$, let $r(M)$, $\nu(M)$ and $c(M)$
denote its rank, nullity and number of coloops. A complex is \emph{stable} if
its Stanley--Reisner ideal is squarefree stable; shifted complexes become
stable after reversing the order of their vertices.

\begin{theorem}[Theorems~\ref{thm:top}, \ref{thm:matroid}
and~\ref{thm:stable}]\label{main:classes}
Let $K$ be a nonvoid complex and let $t\ge0$.
\begin{enumerate}[label=\textup{(\arabic*)}]
\item If $\Ht_{\dim K}(K;\kk)\ne0$ and $t<\hgt I_K$, then
$C(\T_t(K))=L_\kk(\T_t(K))=\dim K+1+t=L_\kk(K)+t$.
\item If $K$ is the independence complex of a matroid $M$ and $t<\nu(M)$,
then $C(\T_t(K))=L_\kk(\T_t(K))=r(M)-c(M)+t$.
\item If $K$ is stable, then
$C(\T_t(K))=L_\kk(\T_t(K))=h(\T_t(K))\le C(K)+t$.
\end{enumerate}
\end{theorem}

Part~(2) is known \cite{KL24,MN24}, \cite[Theorem~30]{BC22}; we include it because it follows directly from part~(1).%
For shifted complexes, we determine the minimal nonfaces of $\T_t(K)$
explicitly. This yields the minimal generators, graded Betti numbers and
regularity of all squarefree symbolic powers of squarefree strongly stable
ideals (Theorem~\ref{thm:shiftexact} and Corollary~\ref{cor:stable}).

Finally, we consider the case $t=1$ for arbitrary $d$, highlighted by Kim
and Lew \cite{KL23}. Since $\eta(r,2)=\lfloor(r+2)^2/4\rfloor$, the
conjectured bound is $\lfloor(d+3)^2/4\rfloor-1$, roughly $d^2/4$, whereas
their general bound is $d^2+2d$. We reduce the conjectured bound for
$d$-collapsible complexes to a single inequality for intersections of links
(Question~\ref{q:links}). A positive answer would show that, for every
$d$-collapsible complex $K$,
$\reg(S/\sym{I_K}{2})\le\left\lfloor\frac{(d+3)^2}{4}\right\rfloor-1$
$(\text{Theorem~\ref{thm:t1conditional}}).$ The answer is positive for $e\le1$, so this holds unconditionally for $d\le2$ and recovers \cite[Theorem~1.6]{KL23}; the first open case is $d=3$.

Table~\ref{tab:summary} summarizes these results.

\begin{table}[ht]
\centering
\caption{Results on Conjectures~\ref{conj:leray} and~\ref{conj:coll}. Here
$C$, $L_\kk$ and $h$ are the invariants of $\T_t(K)$.}
\label{tab:summary}
\small
\renewcommand{\arraystretch}{1.0}
\setlength{\tabcolsep}{6pt}
\begin{tabularx}{\linewidth}{@{}>{\raggedright\arraybackslash}X l l l@{}}
\toprule
Complex $K$ & $t$ & Result & Ref.\\
\midrule
\multicolumn{4}{@{}l}{\textit{Both conjectures}}\\
$m\le d+1$ facets & all & $C\le m(t+1)-1$ & \ref{cor:facets}\\
$\Ht_{\dim K}(K;\kk)\ne0$ & $<\hgt I_K$ & $C=L_\kk=L_\kk(K)+t$ & \ref{thm:top}\\
$\Ind(M)$, $M$ a matroid (known) & $<\nu(M)$ & $C=L_\kk=r(M)-c(M)+t$ & \ref{thm:matroid}\\
stable (e.g.\ shifted, after reversing the order) & all & $C=L_\kk=h\le C(K)+t$ & \ref{thm:stable}\\
\midrule
\multicolumn{4}{@{}l}{\textit{Collapsibility conjecture for $d=1$}}\\
$\Cl(H)$, $H$ a subtree graph of a spider & all & $C\le 2t+1$ & \ref{thm:spider}\\
$1$-collapsible & $\le3$ & $C\le 2t+1$ & \ref{thm:smallt}\\
$1$-collapsible, $\dim K\le4$ & all & $C\le 2t+1$ & \ref{thm:smallt}\\
\midrule
\multicolumn{4}{@{}l}{\textit{Leray conjecture}}\\
flag & $\le3$ & $L_\kk\le L_\kk(K)+2t$ & \ref{cor:flag}\\
$\Ind(G)$, $G$ bipartite or chordal & all & $L_\kk\le L_\kk(K)+2t$ & \ref{cor:flag}\\
$\Ind(G)$, odd girth of $G\ge 2t+1$ & $\ge1$ & $L_\kk\le L_\kk(K)+2t$ & \ref{cor:flag}\\
$K_1*K_2$, Conj.~\ref{conj:leray} for $K_1,K_2$ & all & Conj.~\ref{conj:leray} holds & \ref{cor:joins}\\
\midrule
\multicolumn{4}{@{}l}{\textit{Conditional on a positive answer to Question~\ref{q:links} (unconditional for $d\le2$)}}\\
$d$-collapsible & $1$ & $L_\kk\le\lfloor(d+3)^2/4\rfloor-1$ & \ref{thm:t1conditional}\\
\bottomrule
\end{tabularx}
\end{table}

The paper is organized as follows. Section~\ref{sec:prelim} collects preliminaries, and Section~\ref{sec:sym}
develops the algebraic viewpoint and treats flag complexes.
Sections~\ref{sec:skeleta} and~\ref{sec:config} treat the case $d=1$,
Section~\ref{sec:tone} the case $t=1$, Section~\ref{sec:joins} joins,
Section~\ref{sec:top} top-dimensional homology and matroids, and
Section~\ref{sec:shifted} stable and shifted complexes.

\section{Preliminaries}\label{sec:prelim}

\subsection{Complexes and collapsibility}

A \emph{simplicial complex} on a finite set $W$ is a family
$X\subseteq 2^W$ closed under taking subsets. We allow the
\emph{void complex} $\emptyset$, which has no faces, and
distinguish it from $\{\emptyset\}$. The vertex set of $X$ is
$V(X)=\{v\in W:\{v\}\in X\}$, and the elements of
$W\setminus V(X)$ are its \emph{nonvertices}. For $v\in W$,
the \emph{deletion} and \emph{link} of $v$ in $X$ are
\[
X\setminus v=\{\sigma\in X:v\notin\sigma\},
\qquad
\lk(X,v)=\{\sigma\subseteq W\setminus\{v\}:
\sigma\cup\{v\}\in X\}.
\]
In particular, $\lk(X,v)=\emptyset$ if and only if
$v\notin V(X)$. More generally, for $\sigma\subseteq W$ we put $\lk(X,\sigma)=\{\rho\subseteq W\setminus\sigma:\rho\cup\sigma\in X\}$; it is void if $\sigma\notin X$, and $\lk(X,\{v\})=\lk(X,v)$.
 Deletion and link commute with unions
and intersections. For $U\subseteq W$, the
\emph{induced subcomplex} of $X$ on $U$ is
$X[U]=\{\sigma\in X:\sigma\subseteq U\}$.
We write $2^U$ for the full simplex on $U$.
The dimension of a nonvoid complex $X$ is
$\dim X=\max\{|\sigma|-1:\sigma\in X\}$.
A face of $X$ is \emph{free} if it is contained in
a unique maximal face of $X$. An
\emph{elementary $d$-collapse} removes a free face
$\sigma$ with $|\sigma|\le d$, together with all
faces containing $\sigma$. The complex $X$ is
\emph{$d$-collapsible} if it can be reduced to the
void complex by elementary $d$-collapses. The
\emph{collapsibility number} $C(X)$ is the least
integer $d\ge0$ for which $X$ is $d$-collapsible.
Thus $C(\emptyset)=0$, and every full simplex has
collapsibility number $0$, since its empty face
is free. Every nonvoid complex $X$ satisfies
\begin{equation}\label{eq:dimC}
C(X)\le\dim X+1,
\end{equation}
since its maximal faces can be removed successively.
Fix a field $\kk$. The \emph{Leray number}
$L_\kk(X)$ is the least integer $d\ge0$ such
that $\Ht_i(X[U];\kk)=0$ for every $U\subseteq W$
and every $i\ge d$. We call $X$
\emph{$d$-Leray} over $\kk$ if $L_\kk(X)\le d$.
A \emph{missing face} of $X$ is a set
$\sigma\subseteq W$ such that $\sigma\notin X$
but every proper subset of $\sigma$ belongs to $X$.
The \emph{Helly dimension} of $X$ is
\[
h(X)=\max\bigl(\{0\}\cup
\{|\sigma|-1:\sigma\text{ is a missing face of }X\}\bigr).
\]

The invariants $C(X)$, $L_\kk(X)$, and $h(X)$
are independent of the ambient set $W$.
Indeed, adjoining a nonvertex does not affect
$C(X)$ or $L_\kk(X)$ and creates at most one
new missing face, of cardinality one.
By contrast, tolerance complexes may depend on $W$.

\begin{lemma}[{\cite{Weg75};
see \cite[Section~1]{KL23}}]\label{lem:hLC}
For every simplicial complex $X$,
$h(X)\le L_\kk(X)\le C(X).$
\end{lemma}

We also use the fact that neither $L_\kk$ nor
$C$ increases upon passing to an induced
subcomplex; for $C$, see \cite[Lemma~2.6]{KL23}.
Moreover, for $v\in W$ and $U\subseteq W\setminus\{v\}$ we have
$X[U\cup\{v\}]=X[U]\cup\bigl(v*\lk(X,v)[U]\bigr)$, with intersection
$\lk(X,v)[U]$. The Mayer--Vietoris sequence therefore shows that
$L_\kk(\lk(X,v))\le L_\kk(X)$ and
\begin{equation}\label{eq:LerayDL}
L_\kk(X)\le\max\{L_\kk(X\setminus v),\,L_\kk(\lk(X,v))+1\}.
\end{equation}
Indeed, if $\lk(X,v)[U]$ is nonvoid, then for $i\ge0$ the exact sequences $\Ht_{i+1}(X[U\cup\{v\}];\kk)\to\Ht_i(\lk(X,v)[U];\kk)\to\Ht_i(X[U];\kk)$ and $\Ht_i(X[U];\kk)\to\Ht_i(X[U\cup\{v\}];\kk)\to\Ht_{i-1}(\lk(X,v)[U];\kk)$ give the two inequalities; if $\lk(X,v)[U]$ is void, then $X[U\cup\{v\}]=X[U]$.
The following inequality is our main tool
for bounding collapsibility numbers.

\begin{lemma}[{Tancer \cite[Proposition~1.2]{Tan11};
see \cite[Lemma~3.1]{KL22}}]\label{lem:tancer}
Let $X$ be a simplicial complex and let
$v\in V(X)$. Then
$C(X)\le\max\bigl\{
C(X\setminus v),\ C(\lk(X,v))+1
\bigr\}.
$
In particular, if $d\ge1$, $X\setminus v$ is
$d$-collapsible, and $\lk(X,v)$ is
$(d-1)$-collapsible, then $X$ is
$d$-collapsible.
\end{lemma}

\subsection{Tolerance complexes}

For complexes $X,Y$ on disjoint sets, their \emph{join} is
$X*Y=\{\sigma\cup\rho:\sigma\in X,\ \rho\in Y\}$.
For $\tau\subseteq W$ and $r\in\mathbb Z$, the
\emph{tolerance ball} of radius $r$ centered at $\tau$ is
\[
\Sk_r(\tau)=\{\sigma\subseteq W:|\sigma\setminus\tau|\le r\}
=2^\tau*\{\rho\subseteq W\setminus\tau:|\rho|\le r\}.
\]
It is void for $r<0$, and $\Sk_0(\tau)=2^\tau$.
For a complex $K$ on $W$ and $s\in\mathbb Z$, define
\[
\T_s(K)=\{\sigma\subseteq W:|\sigma\setminus f|\le s
\text{ for some }f\in K,\ f\subseteq\sigma\}.
\]
Thus $\T_s(K)$ is void for $s<0$, and $\T_0(K)=K$.
Writing $\mathcal F(K)$ for the set of facets of $K$,
we have
\begin{equation}\label{eq:Tballs}
\T_s(K)=\bigcup_{F\in\mathcal F(K)}\Sk_s(F).
\end{equation}
Indeed, extending $f$ to a facet $F$ does not increase
$|\sigma\setminus f|$. Conversely, if
$|\sigma\setminus F|\le s$, then
$f=\sigma\cap F\in K$ satisfies
$|\sigma\setminus f|\le s$.
In particular, $\Sk_r(\tau)=\T_r(2^\tau)$.
Moreover,
$\T_s(\T_t(K))=\T_{s+t}(K)$ for all $s,t\ge0$,
as follows by partitioning $\sigma\setminus f$
into sets of sizes at most $t$ and $s$: if $f\subseteq\sigma$ and $|\sigma\setminus f|\le s+t$, choose $g$ with $f\subseteq g\subseteq\sigma$, $|g\setminus f|\le t$ and $|\sigma\setminus g|\le s$; then $g\in\T_t(K)$.

\begin{lemma}\label{lem:dom}
Let $\tau,\tau'\subseteq W$ and $r,r'\in\mathbb Z$.
If $|\tau\setminus\tau'|+r\le r'$, then
$\Sk_r(\tau)\subseteq\Sk_{r'}(\tau')$.
\end{lemma}

\begin{proof}
For every $\sigma\subseteq W$,
$|\sigma\setminus\tau'|
\le|\sigma\setminus\tau|+|\tau\setminus\tau'|.$
The assertion follows.
\end{proof}

\begin{lemma}\label{lem:Sdl}
Let $\tau\subseteq W$, $v\in W$, and $r\in\mathbb Z$.
As complexes on $W\setminus\{v\}$:
\begin{enumerate}[label=\textup{(\arabic*)}]
\item $\Sk_r(\tau)\setminus v
=\Sk_r(\tau\setminus\{v\})$;
\item
$\lk(\Sk_r(\tau),v)=
\begin{cases}
\Sk_r(\tau\setminus\{v\}),&\text{ if }v\in\tau,\\
\Sk_{r-1}(\tau),&\text{ if }v\notin\tau.
\end{cases}
$
\end{enumerate}
\end{lemma}

\begin{proof}
For $\sigma\subseteq W\setminus\{v\}$, we have
$|\sigma\setminus\tau|
=|\sigma\setminus(\tau\setminus\{v\})|$.
Moreover,
$|(\sigma\cup\{v\})\setminus\tau|$
equals $|\sigma\setminus\tau|$ if $v\in\tau$,
and $|\sigma\setminus\tau|+1$ otherwise.
Both identities follow.
\end{proof}

\begin{lemma}[{cf.\ \cite[Section~4]{KL23}}]\label{lem:Tdl}
Let $K$ be a complex on $W$, $v\in W$, and $s\in\mathbb Z$. With all
tolerance complexes on the right formed on $W\setminus\{v\}$,
\[
\T_s(K)\setminus v=\T_s(K\setminus v),
~
\lk(\T_s(K),v)=\T_{s-1}(K\setminus v)\cup\T_s(\lk(K,v)).
\]
\end{lemma}

\begin{proof}
The first identity holds because $v\notin\sigma$ and $f\subseteq\sigma$
imply $v\notin f$. For the second, let $\sigma\subseteq W\setminus\{v\}$
and suppose that $\sigma\cup\{v\}\in\T_s(K)$, witnessed by $f\in K$ with
$f\subseteq\sigma\cup\{v\}$ and $|(\sigma\cup\{v\})\setminus f|\le s$.
If $v\notin f$, then $f\in K\setminus v$ and $|\sigma\setminus f|\le s-1$.
If $v\in f$, then $f\setminus\{v\}\in\lk(K,v)$ and
$|\sigma\setminus(f\setminus\{v\})|\le s$. This proves one inclusion.
Conversely, a witnessing face of $\sigma$ in $K\setminus v$ for
$\T_{s-1}(K\setminus v)$ also witnesses $\sigma\cup\{v\}\in\T_s(K)$, and so
does the union of $\{v\}$ with a witnessing face of $\sigma$ in $\lk(K,v)$
for $\T_s(\lk(K,v))$.
\end{proof}

\begin{lemma}[{cf.\ \cite[Proof of Theorem~1.5]{KL23}}]
\label{lem:Tind}
Let $K$ be a complex on $W$, $s\in\mathbb Z$,
and $U\subseteq W$. Then
$\T_s(K)[U]=\T_s(K[U])$, where the right-hand
side is formed on $U$.
\end{lemma}

\begin{proof}
If $f\subseteq\sigma\subseteq U$ and $f\in K$, then $f\in K[U]$.
\end{proof}

\subsection{Erd\H{o}s--Gallai numbers}

Let $\cH$ be a finite hypergraph with nonempty
edges. Its \emph{transversal number} $\tau(\cH)$
is the minimum cardinality of a set meeting
every edge. A hypergraph without isolated vertices
is \emph{$k$-critical} if $\tau(\cH)=k$ and
$\tau(\cH-e)<k$ for every edge $e$.
The \emph{Erd\H{o}s--Gallai number} $\eta(r,k)$
is the maximum number of vertices of an
$r$-uniform $k$-critical hypergraph.
Clearly $\eta(1,k)=k$, and Erd\H{o}s and Gallai
\cite{EG61} proved that
\[
\eta(2,k)=2k,~
\eta(r,2)=\left\lfloor\frac{(r+2)^2}{4}\right\rfloor.
\]
Tuza \cite{Tuz85} proved that, for $r\ge2$ and
$k\ge1$,
$\eta(r,k)<
\binom{r+k-1}{r-1}+\binom{r+k-2}{r-1}.$
In particular, $\eta(r,k)$ is finite for all
$r,k\ge1$.

\begin{lemma}\label{lem:eta}
Let $r,k\ge1$.
\begin{enumerate}[label=\textup{(\arabic*)}]
\item $rk\le\eta(r,k)\le\eta(r+1,k)$.
\item Every $k$-critical hypergraph whose edges
have cardinality at most $r$ has at most
$\eta(r,k)$ vertices.
\end{enumerate}
\end{lemma}

\begin{proof}
We first observe that adjoining new vertices to
the edges of a hypergraph, with each new vertex
belonging to exactly one edge, preserves its
transversal number. Indeed, every transversal
of the original hypergraph remains a transversal,
while each new vertex in a transversal of the
enlarged hypergraph can be replaced by a vertex
of its original edge. The same argument applies
after deleting any edge. Thus this operation
preserves $k$-criticality. 
(1) A hypergraph consisting of $k$ pairwise
disjoint edges of size $r$ is $k$-critical
and has $rk$ vertices. The second inequality
follows by adjoining one new vertex to each
edge of an $r$-uniform $k$-critical hypergraph. 
(2) In a $k$-critical hypergraph no edge contains another, since deleting the larger edge would not lower the transversal number; hence the edges below remain distinct. Adjoin distinct new vertices to the edges
until each edge has cardinality $r$.
The resulting hypergraph is $r$-uniform and
$k$-critical, and has at least as many vertices
as the original hypergraph.
\end{proof}

\subsection{Stanley--Reisner rings and squarefree symbolic powers}

Let $K$ be a simplicial complex on $W$, let
$S=\kk[x_w:w\in W]$, and write
$x^\sigma=\prod_{w\in\sigma}x_w$ for
$\sigma\subseteq W$. The
\emph{Stanley--Reisner ideal} of $K$ is
$I_K=(x^\sigma:\sigma\subseteq W,\ \sigma\notin K).$
Its minimal monomial generators correspond to
the missing faces of $K$. For $A\subseteq W$,
let $P_A=(x_w:w\in A)$. If $K$ is nonvoid, then
\[
I_K=\bigcap_{F\in\mathcal F(K)}Q_F,
\qquad Q_F=P_{W\setminus F},
\]
is the irredundant decomposition of $I_K$
into minimal primes. Consequently,
$\dim S/I_K=\dim K+1$ and
$\hgt I_K=|W|-\dim K-1$.

We recall some homological invariants of
finitely generated graded $S$-modules.
Let $M$ be such a module, with minimal
graded free resolution
\[
0\longrightarrow F_p\longrightarrow
\cdots\longrightarrow F_1
\longrightarrow F_0
\longrightarrow M\longrightarrow 0,
\]
where
$F_i=\bigoplus_{j\in\mathbb Z}
S(-j)^{\beta_{i,j}(M)}.$
The integers $\beta_{i,j}(M)$ are the
\emph{graded Betti numbers} of $M$.
The \emph{projective dimension} and
\emph{Castelnuovo--Mumford regularity}
of a nonzero module $M$ are defined by
$\pd_S M
=\max\{i:\beta_{i,j}(M)\ne0
\text{ for some }j\}$
and
$\reg M
=\max\{j-i:\beta_{i,j}(M)\ne0\},$
respectively. 
For a multigraded $S$-module $M$,
its \emph{multigraded Betti numbers} are
given by
$\beta_{i,\mathbf a}(M)
=\dim_\kk\Tor_i^S(M,\kk)_{\mathbf a},$
$\mathbf a\in\mathbb N^W.$
Identifying a subset $\sigma\subseteq W$
with its squarefree multidegree,
Hochster's formula \cite{Hoc77} gives
\[
\beta_{i,\sigma}(S/I_K)
=\dim_\kk\Ht_{|\sigma|-i-1}(K[\sigma];\kk),
~ \sigma\subseteq W,
\]
with the convention
$\Ht_{-1}(\{\emptyset\};\kk)=\kk$.
Moreover, these Betti numbers vanish
in nonsquarefree multidegrees, and
$\beta_{i,j}(S/I_K)
=\sum_{\substack{\sigma\subseteq W\\
|\sigma|=j}}
\beta_{i,\sigma}(S/I_K).$
Hence $\reg(S/I_K)=L_\kk(K)$ for every
nonvoid $K$; see \cite[Section~1]{KM06}.
If moreover $K\ne2^W$, then
$\reg I_K=L_\kk(K)+1$, while $h(K)+1$
is the largest degree of a minimal monomial
generator of $I_K$.

Let $I=\bigcap_j P_j$ be a squarefree monomial ideal,
where the $P_j$ are its minimal primes. For $k\ge1$,
its $k$-th symbolic power is
$I^{(k)}=\bigcap_j P_j^{\,k}$
\cite[Proposition~1.4.4]{HH11}.
For a monomial ideal $J$, let $J_\sqf$ denote the ideal
generated by its squarefree monomials. The
\emph{$k$-th squarefree symbolic power} of $I$ is
$\sym{I}{k}=(I^{(k)})_\sqf;$
see \cite{SF23,CDRS}.
For $A\subseteq W$, put
$P_A^{[k]}=(P_A^k)_\sqf$.
This ideal is generated by the squarefree monomials
of degree $k$ in the variables of $P_A$, and is zero
if $k>|A|$. By \cite[Lemma~2.3]{CDRS},
$\sym{I}{k}=\bigcap_j P_j^{[k]}.$
In particular, the zero ideal, whose only minimal
prime is $P_\emptyset=0$, satisfies
$\sym{0}{k}=0$.

\section{Squarefree symbolic powers}\label{sec:sym}

In this section, we identify the Stanley--Reisner ideals of
tolerance complexes with squarefree symbolic powers, describe
their minimal generators in terms of transversal numbers, and
reformulate Conjecture~\ref{conj:leray} in terms of regularity.
We then use a comparison with symbolic powers to establish
the conjecture for certain flag complexes. Unless otherwise
stated, $K$ is a nonvoid simplicial complex on a finite set
$W$, and $I=I_K$.

\begin{proposition}\label{prop:sym}
For every integer $t\ge0$,
$I_{\T_t(K)}=\sym{I_K}{t+1}
=\bigcap_{F\in\mathcal F(K)}Q_F^{[t+1]}.
$
\end{proposition}

\begin{proof}
The second equality follows from
\cite[Lemma~2.3]{CDRS}. For the first, let
$\sigma\subseteq W$. By \eqref{eq:Tballs},
\[
\begin{aligned}
\sigma\notin\T_t(K)
&\iff |\sigma\setminus F|\ge t+1
\text{ for every }F\in\mathcal F(K)
\iff x^\sigma\in
\bigcap_{F\in\mathcal F(K)}Q_F^{[t+1]}.
\end{aligned}
\]
Since both ideals are squarefree monomial ideals,
they coincide.
\end{proof}

Mantero and Nguyen \cite[Proposition~3.19]{MN24}
give an equivalent description of squarefree
symbolic powers in terms of cover ideals. 
To describe the minimal generators, let $\cH$
be the hypergraph on $W$ whose edges are the
missing faces of $K$. For $\sigma\subseteq W$,
let $\cH[\sigma]$ denote the hypergraph whose
edges are those of $\cH$ contained in $\sigma$.

\begin{proposition}\label{prop:gens}
For every integer $k\ge1$ and every
$\sigma\subseteq W$,
\[
x^\sigma\in\sym{I}{k}
\iff \sigma\notin\T_{k-1}(K)
\iff \tau(\cH[\sigma])\ge k.
\]
Consequently, the minimal monomial generators
of $\sym{I}{k}$ are the monomials $x^\sigma$
for which $\sigma$ is inclusion-minimal
subject to $\tau(\cH[\sigma])\ge k$.
For every such $\sigma$,
$\tau(\cH[\sigma])=k$.
\end{proposition}

\begin{proof}
The first equivalence follows from
Proposition~\ref{prop:sym}. For the second,
a subset $f\subseteq\sigma$ belongs to $K$
if and only if $\sigma\setminus f$ is a
transversal of $\cH[\sigma]$. Thus
$\tau(\cH[\sigma])
=\min\{|\sigma\setminus f|:
f\in K,\ f\subseteq\sigma\},
$
where the transversal number of a hypergraph
with no edges is zero. The description of
the minimal generators follows. 
Finally, deleting a vertex decreases the
transversal number by at most one, since adding the vertex to a transversal of the smaller hypergraph gives a transversal of the larger one; note that $\sigma\ne\emptyset$ because $k\ge1$. Hence
inclusion-minimality forces
$\tau(\cH[\sigma])=k$.
\end{proof}

\begin{proposition}\label{prop:height}
Let $\omega=\dim K+1$, so that
$\hgt I=|W|-\omega$. If $0\le t<\hgt I$, then
\[
\dim\T_t(K)=\dim K+t,~
\hgt\sym{I}{t+1}=\hgt I-t,
\]
and
$\dim S/\sym{I}{t+1}=\dim S/I+t$.
If $t\ge\hgt I$, then
$\T_t(K)=2^W$ and $\sym{I}{t+1}=0$.
\end{proposition}

\begin{proof}
Let $F$ be a facet of $K$ with
$|F|=\omega$. Every face of $\T_t(K)$
has cardinality at most $\omega+t$.
If $t<\hgt I=|W\setminus F|$, choose
$A\subseteq W\setminus F$ with $|A|=t$.
Then $F\cup A\in\T_t(K)$, so
$\dim\T_t(K)=\dim K+t$. The remaining
identities follow from
Proposition~\ref{prop:sym} and the
Stanley--Reisner correspondence. 
If $t\ge\hgt I$, then
$|\sigma\setminus F|\le|W\setminus F|\le t$
for every $\sigma\subseteq W$. Hence
$\T_t(K)=2^W$ and $\sym{I}{t+1}=0$.
\end{proof}

In particular, $\sym{I}{k}\ne0$ if and only if
$k\le\hgt I$; see \cite[Proposition~2.2]{SF23}. 
By Proposition~\ref{prop:sym} and Hochster's
formula,
$L_\kk(\T_t(K))=\reg(S/\sym{I_K}{t+1})$.
Thus Conjecture~\ref{conj:leray} is equivalent
to the following algebraic statement.

\begin{conjecture}[Algebraic form of
Conjecture~\ref{conj:leray}]\label{conj:reg}
Let $I\subsetneq S$ be a squarefree monomial
ideal and let $d\ge0$. If $\reg(S/I)\le d$,
then, for every integer $t\ge0$,
\[
\reg(S/\sym{I}{t+1})
\le\eta(d+1,t+1)-1.
\]
\end{conjecture}

By Lemma~\ref{lem:hLC},
Conjecture~\ref{conj:coll} implies the
conclusion of Conjecture~\ref{conj:reg}
whenever $I=I_K$ and $K$ is $d$-collapsible.
In algebraic terms, Kim and Lew
\cite[Theorems~1.5 and~1.6]{KL23} proved that
$\reg(S/\sym{I_K}{t+1})\le\hKL(t,d)$
whenever $K$ is $d$-collapsible, and that
$\reg(S/\sym{I_K}{2})\le5$ whenever
$K$ is $2$-collapsible. Both bounds hold
over any field.
For edge ideals, Seyed Fakhari conjectured that
$\reg\sym{I(G)}{s}\le\operatorname{match}(G)+s$
whenever $\sym{I(G)}{s}\ne0$, where
$\operatorname{match}(G)$ is the matching number
of $G$ \cite[Conjecture~1.1]{SF23}. Equivalently,
$L_\kk(\T_t(\Ind(G)))\le\operatorname{match}(G)+t$
whenever $\T_t(\Ind(G))\ne2^{V(G)}$. Since $\operatorname{match}(G)$ can exceed $L_\kk(\Ind(G))=\reg I(G)-1$, this would not by itself imply Conjecture~\ref{conj:leray} for flag complexes.

\begin{remark}\label{rem:dimbound}
Every face of $\T_t(K)$ has cardinality
at most $\dim K+1+t$. Hence, by
\eqref{eq:dimC},
$C(\T_t(K))\le\dim K+1+t$, and consequently
$\reg(S/\sym{I}{t+1})\le\dim S/I+t$.
\end{remark}

We now compare squarefree symbolic powers with symbolic powers. For a
monomial ideal $J$, the ideal $J_\sqf$ is the restriction of $J$ to the
multidegree $x^W$, so the restriction lemma \cite[Lemma~4.4]{HHZ04}
applies. Erey, Herzog, Hibi, and Saeedi Madani \cite[Lemma~1.2]{EHHS22}
used it to show that upper bounds on the regularity of ordinary powers pass
to squarefree powers. The same holds for symbolic powers.

\begin{proposition}\label{prop:compare}
Let $I$ be a squarefree monomial ideal and let $k\ge1$ with
$\sym{I}{k}\ne0$. Then
$\beta_{i,\sigma}(\sym{I}{k})=\beta_{i,\sigma}(I^{(k)})$ for every $i$ and
every $\sigma\subseteq W$. Consequently, $\reg\sym{I}{k}\le\reg I^{(k)}$.
\end{proposition}

\begin{proof}
Since $\sym{I}{k}$ is generated by the minimal generators of $I^{(k)}$
dividing $x^W$, this follows from the restriction lemma
\cite[Lemma~4.4]{HHZ04}; the regularity inequality is also a special case
of \cite[Proposition~2.3]{SF23}.
\end{proof}

For a graph $G$, let $I(G)$ denote its edge ideal and $\Ind(G)$ its
independence complex. Thus $I(G)=I_{\Ind(G)}$ on $V(G)$, and every flag
complex is the independence complex of a graph.

\begin{theorem}\label{thm:flag}
Let $G$ be a graph with at least one edge, and let $t\ge0$. Suppose that
\begin{equation}\label{eq:symbound}
\reg I(G)^{(s)}\le\reg I(G)+2s-2\qquad(1\le s\le t+1).
\end{equation}
Let $W\supseteq V(G)$ be finite, and regard $K=\Ind(G)$ as a complex on $W$
for which the elements of $W\setminus V(G)$ are nonvertices. Then
$L_\kk(\T_s(K))\le L_\kk(K)+2s$ for $0\le s\le t$. Consequently,
Conjecture~\ref{conj:leray} holds for $K$ and this value of $t$.
\end{theorem}

\begin{proof}
First suppose $W=V(G)$. The assertion is immediate for $s=0$. Let
$1\le s\le t$. If $\sym{I(G)}{s+1}=0$, then $\T_s(K)=2^W$ and the bound
holds. Otherwise, Propositions~\ref{prop:sym} and~\ref{prop:compare},
together with~\eqref{eq:symbound}, give
\[
L_\kk(\T_s(K))=\reg\sym{I(G)}{s+1}-1
\le\reg I(G)^{(s+1)}-1
\le\reg I(G)+2s-1
=L_\kk(K)+2s,
\]
where the last equality holds because $K\ne2^W$, so
$\reg I(G)=L_\kk(K)+1$.

We extend the bound to $W\supseteq V(G)$ by induction on
$|W\setminus V(G)|$, simultaneously for $0\le s\le t$. The case $s=0$
follows from the invariance of the Leray number under adjoining
nonvertices. Suppose $s\ge1$, choose $z\in W\setminus V(G)$, and put
$U=W\setminus\{z\}$ and $K'=K[U]$; note that $L_\kk(K')=L_\kk(K)$. Since $z$
is a nonvertex of $K$, Lemma~\ref{lem:Tdl} gives
$\T_s(K)\setminus z=\T_s(K')$, $\lk(\T_s(K),z)=\T_{s-1}(K').$
By~\eqref{eq:LerayDL} applied to $X=\T_s(K)$ and $v=z$, this yields
$L_\kk(\T_s(K))\le\max\{L_\kk(\T_s(K')),\,L_\kk(\T_{s-1}(K'))+1\}.$
By induction, the right-hand side is at most $L_\kk(K)+2s$.
Finally, suppose that $K$ is $d$-Leray. Since $G$ has an edge, $K$ has a
missing face of cardinality $2$, so $d\ge L_\kk(K)\ge h(K)\ge1$ by
Lemma~\ref{lem:hLC}. By Lemma~\ref{lem:eta}(1),
$L_\kk(\T_t(K))\le d+2t\le(d+1)(t+1)-1\le\eta(d+1,t+1)-1.$
Thus Conjecture~\ref{conj:leray} holds for $K$ and $t$.
\end{proof}

\begin{corollary}\label{cor:flag}
Let $G$ be a graph with at least one edge,
and let $K=\Ind(G)$ on a finite set
$W\supseteq V(G)$, for which the elements of $W\setminus V(G)$ are nonvertices. Then
Conjecture~\ref{conj:leray} holds for $K$,
with $L_\kk(\T_t(K))\le L_\kk(K)+2t$,
in each of the following cases:
\begin{enumerate}[label=\textup{(\arabic*)}]
\item $t\le3$, for every graph $G$;
\item all $t$, if $G$ is bipartite;
\item all $t$, if $G$ is chordal;
\item $t\ge1$, if $G$ has no odd cycle of
length at most $2t-1$.
\end{enumerate}
% Moreover, in each case,
% $\reg\sym{I(G)}{s}\le\reg I(G)+2s-2$
% for every $1\le s\le t+1$ with
% $\sym{I(G)}{s}\ne0$.
\end{corollary}

\begin{proof}
Inequality~\eqref{eq:symbound} holds in case (1) by
\cite[Theorem~1.3]{MNPTV22}. In cases (2) and (4) with $t\ge1$, it holds by
\cite[Corollary~3.11]{SF19b} applied with $k=t$: if $G$ has no odd cycle of
length at most $2t-1$, which is automatic when $G$ is bipartite, then
$\reg I(G)^{(s)}\le\reg I(G)+2s-2$ for all $s\le t+1$. For $t=0$,
inequality~\eqref{eq:symbound} is trivial.
For case (3), Seyed Fakhari \cite{SF22} proved that
$\reg I(G)^{(s)}=2s+\operatorname{indmatch}(G)-1$ for every $s\ge1$, where
$\operatorname{indmatch}(G)$ is the induced matching number of $G$. Taking
$s=1$ gives $\reg I(G)=\operatorname{indmatch}(G)+1$, so
inequality~\eqref{eq:symbound} holds with equality.
The conclusions now follow from Theorem~\ref{thm:flag}.
\end{proof}

Thus Corollary~\ref{cor:flag} is a direct translation, via Theorem~\ref{thm:flag}, of known regularity bounds for symbolic powers of edge ideals.

If $G$ has no edges, then $K=2^{V(G)}$ has a unique facet, and
Corollary~\ref{cor:facets} applies. Moreover, in this case
$\T_t(K)=\mathrm B_t(V(G))
=2^{V(G)}*\{\rho\subseteq W\setminus V(G):|\rho|\le t\}$, whose Leray
number is $\min\{t,|W\setminus V(G)|\}$, since it is the join of a simplex with the $(t-1)$-skeleton of the simplex on $W\setminus V(G)$. This is at most $t\le L_\kk(K)+2t$. Thus Conjecture~\ref{conj:leray}
holds for all flag complexes when $t\le3$, and the bound
$L_\kk(\T_t(K))\le L_\kk(K)+2t$ holds for every flag complex $K$ in the
cases of Corollary~\ref{cor:flag}.
For $d=1$, the conjecture also follows from \cite[Theorem~1.5]{KL23},
since flag $1$-Leray complexes are $1$-collapsible by
Lemma~\ref{lem:chordal}. For $d\ge2$, that theorem gives only the bound
$h_{KL}(t,d)$, which is in general larger than $\eta(d+1,t+1)-1$ (for
instance, $h_{KL}(1,2)=8>5$; the case $d=2$, $t=1$ is covered by
\cite[Theorem~1.6]{KL23}). It also applies only to $d$-collapsible
complexes, whereas Corollary~\ref{cor:flag} gives the conjectured bound for
every flag $d$-Leray complex when $t\le3$.

\section{Unions of tolerance balls}\label{sec:skeleta}

In this section, we bound the collapsibility number of
unions of tolerance balls using deletion--link induction
(Lemma~\ref{lem:tancer}). The corresponding Leray bounds
follow in certain cases from the union theorem of
Kalai and Meshulam \cite{KM06}, but the collapsibility
bounds require a different argument. These results are
the main tools for the case $d=1$ in
Section~\ref{sec:config}.

\begin{definition}
A \emph{term} on $W$ is a pair $(\tau,r)$ with
$\tau\subseteq W$ and $r\in\mathbb Z$.
A \emph{chain} is a nonempty finite list of terms
$\cK=((\tau_1,r_1),\dots,(\tau_q,r_q))$ satisfying
$\tau_1\subseteq\dots\subseteq\tau_q$ and
$r_1\ge\dots\ge r_q$. Its \emph{top} is
$\topp(\cK)=r_1$.
A \emph{system} is a finite list of chains
$\cS=(\cK_1,\dots,\cK_m)$.
For a finite list $\cL$ of terms, and in particular
for a chain or a system, let $Y(\cL)$ denote the
union of the balls $\Sk_r(\tau)$ over its terms.
\end{definition}

Every nonempty sublist of a chain is a chain, with
top at most that of the original chain. Moreover,
if $(\tau,r)$ and $(\tau',r')$ occur in distinct
positions of a list $\cL$ and
$|\tau\setminus\tau'|+r\le r'$, then deleting
$(\tau,r)$ does not change $Y(\cL)$, by
Lemma~\ref{lem:dom}. We call this operation
\emph{clean-up}. The two terms may belong to different chains of a system. Chains that become empty
are discarded. Since the dominating term remains, clean-up never removes every term of a nonempty list.

\begin{theorem}\label{thm:skeleta}
Let $\cS=(\cK_1,\dots,\cK_m)$, $m\ge1$, be a
system on $W$ whose terms all have nonnegative
radius. Then
\[
C(Y(\cS))\le\Psi(\cS):=
\sum_{i=1}^m\bigl(\topp(\cK_i)+1\bigr)-1.
\]
\end{theorem}

\begin{proof}
We induct on $|W|$. Put $Y=Y(\cS)$ and
$d=\Psi(\cS)\ge0$. If some term $(\tau,r)$
satisfies $|W\setminus\tau|\le r$, then
$Y=2^W$ and $C(Y)=0$. This also covers
$W=\emptyset$. 
Otherwise, let $\tau^+$ be the last center
of $\cK_1$ and choose $v\in W\setminus\tau^+$.
By nesting, $v$ belongs to no center of
$\cK_1$.
Replacing every term $(\tau,r)$ by
$(\tau\setminus\{v\},r)$ gives a system
$\cS\setminus v$ on $W\setminus\{v\}$ with
$\Psi(\cS\setminus v)=d$. By
Lemma~\ref{lem:Sdl}(1),
$Y\setminus v=Y(\cS\setminus v)$, so
$C(Y\setminus v)\le d$ by induction.
If $v\notin V(Y)$, then $Y=Y\setminus v$
and we are done. Hence assume $v\in V(Y)$.

Consider a chain $\cK=((\tau_j,r_j))_j$
of $\cS$. Since its centers are nested,
the terms with $v\in\tau_j$ form a final
segment. By Lemma~\ref{lem:Sdl}(2),
$\lk(Y(\cK),v)$ is represented by the list
consisting of $(\tau_j,r_j-1)$ for
$v\notin\tau_j$, followed by
$(\tau_j\setminus\{v\},r_j)$ for
$v\in\tau_j$, after discarding terms
of negative radius. The centers remain
nested, and the radii are nonincreasing
within each segment.

Suppose both segments are nonempty, and
let $(\tau_b,r_b)$ be the first term of
$\cK$ with $v\in\tau_b$. Every earlier
center is contained in $\tau_b\setminus\{v\}$.
Thus, if $j<b$ and $r_j=r_b$, the term
$(\tau_j,r_j-1)$ is dominated by
$(\tau_b\setminus\{v\},r_b)$ and can be
removed by clean-up. These terms form
a terminal part of the first segment.
The remaining radii in that segment
are at least $r_b$, so the resulting
list is a chain.
Applying this construction to each chain
gives a chain $\cK^v$, possibly empty.
After discarding empty chains, we obtain
a system $\cS^v$ on $W\setminus\{v\}$
such that $\lk(Y,v)=Y(\cS^v)$, since links commute with unions.
All its terms have nonnegative radius,
and $\cS^v$ is nonempty because $v\in V(Y)$.

The top of each surviving chain is at most
its original top. For $\cK_1$, all radii
decrease by one, so its contribution
$\topp(\cK_1)+1$ decreases by exactly one;
the chain disappears if its top is $0$.
Consequently,
$0\le\Psi(\cS^v)\le d-1$.
In particular, $d\ge1$, and induction
gives $C(\lk(Y,v))\le d-1$.
Together with $C(Y\setminus v)\le d$,
Lemma~\ref{lem:tancer} yields $C(Y)\le d$.
\end{proof}

\begin{corollary}\label{cor:facets}
If $K$ has $m\ge1$ facets, then
$C(\T_t(K))\le m(t+1)-1$ for every $t\ge0$.
In particular, $C(K)\le m-1$.
If $m\le d+1$, then $\T_t(K)$ is
$(\eta(d+1,t+1)-1)$-collapsible.
Consequently, Conjectures~\ref{conj:leray}
and~\ref{conj:coll} hold for complexes
with at most $d+1$ facets.
\end{corollary}

\begin{proof}
By \eqref{eq:Tballs}, Theorem~\ref{thm:skeleta}
applies to the $m$ single-term chains
$((F,t))$, where $F$ ranges over the
facets of $K$. Taking $t=0$ gives
$C(K)\le m-1$. If $m\le d+1$, then
\[
m(t+1)-1\le(d+1)(t+1)-1
\le\eta(d+1,t+1)-1
\]
by Lemma~\ref{lem:eta}(1). The conclusions
follow from Lemma~\ref{lem:hLC}.
\end{proof}

\begin{example}\label{ex:sharp}
Corollary~\ref{cor:facets} is sharp for
all $m\ge1$ and $t\ge0$. Let
$W=A_1\sqcup\dots\sqcup A_m$, where
$|A_i|=t+1$, and let $K$ have facets
$F_i=W\setminus A_i$. Then
$\sigma\in\T_t(K)$ if and only if
$|\sigma\cap A_i|\le t$ for some $i$,
which is equivalent to $\sigma\ne W$.
Thus $\T_t(K)$ is the boundary of the
simplex on $W$, and
$C(\T_t(K))=L_\kk(\T_t(K))
=h(\T_t(K))=m(t+1)-1.$
The missing faces of $K$ are the minimal
transversals of $\{A_1,\dots,A_m\}$,
so $h(K)=m-1$. Hence
$C(K)=L_\kk(K)=h(K)=m-1$ by
Corollary~\ref{cor:facets} and
Lemma~\ref{lem:hLC}. For $m=2$,
this is the example of Kim and Lew
\cite[Section~1]{KL23} showing that
the bound $2t+1$ is sharp for $d=1$.
\end{example}

\section{The case \texorpdfstring{$d=1$}{d=1}}
\label{sec:config}

In this section, we study Conjecture~\ref{conj:coll}
for $d=1$. We first reduce it to clique complexes of
chordal graphs. We then
introduce a class of complexes closed under deletion
and link at admissible simplicial vertices, yielding
an induction principle controlled by a potential
function. We apply this principle to prove the
conjectured bound for subtree graphs of spiders and for all
chordal graphs when $t\le3$. 
For a graph $H$ on $W$, let $\Cl(H)$ denote its
clique complex. A vertex $u$ of $H$ is
\emph{simplicial} if its neighborhood $N_H(u)$
is a clique. Every noncomplete chordal graph
has two nonadjacent simplicial vertices
\cite{Dir61}.

\begin{lemma}[{see \cite{Weg75,Fro90}}]
\label{lem:chordal}
For a nonvoid complex $K$ on its vertex set $V$,
the following are equivalent:
\textup{(i)} $K$ is $1$-collapsible;
\textup{(ii)} $K$ is $1$-Leray over some field;
\textup{(iii)} $K=\Cl(H)$ for a chordal graph
$H$ on $V$.
\end{lemma}

Here (i) implies (ii) by Lemma~\ref{lem:hLC}; (ii) implies (iii) because a $1$-Leray complex has no missing face of cardinality at least $3$, and an induced cycle of length at least $4$ in its $1$-skeleton would give nonzero $\Ht_1$; and (iii) implies (i) by collapsing at simplicial vertices \cite{Weg75}.

\begin{lemma}\label{lem:ghost}
Let $H$ be a chordal graph and let $t\ge0$.
Suppose that $\T_s(\Cl(H))$ is
$(2s+1)$-collapsible for $0\le s\le t$.
Let $K$ be a complex on a finite set $W$
whose restriction to its vertex set $V(H)$
is $\Cl(H)$ and whose remaining elements
are nonvertices. Then $\T_s(K)$ is
$(2s+1)$-collapsible for $0\le s\le t$.
\end{lemma}

\begin{proof}
We induct on $|W\setminus V(H)|$. The case
$W=V(H)$ holds by assumption. Otherwise,
choose $z\in W\setminus V(H)$ and let
$K'=K[W\setminus\{z\}]$. Since
$\lk(K,z)=\emptyset$, Lemma~\ref{lem:Tdl}
gives
$\T_s(K)\setminus z=\T_s(K'),$ $\lk(\T_s(K),z)=\T_{s-1}(K').$
For $s=0$, the vertex $z$ does not belong
to $\T_0(K)=K$, so the assertion follows
by induction. For $s\ge1$, induction gives
$C(\T_s(K'))\le2s+1$ and
$C(\T_{s-1}(K'))\le2s-1$.
Lemma~\ref{lem:tancer} now yields
$C(\T_s(K))\le2s+1$.
\end{proof}

By Lemmas~\ref{lem:chordal} and~\ref{lem:ghost},
the case $d=1$ of Conjecture~\ref{conj:coll}
is equivalent to
$C(\T_t(\Cl(H)))\le2t+1$ for every chordal
graph $H$ and every $t\ge0$.

\begin{lemma}\label{lem:simplicial}
Let $u$ be a simplicial vertex of $H$,
let $K=\Cl(H)$, put $\rho=N_H(u)$,
and let $s\ge0$. As complexes on
$W\setminus\{u\}$,
\[
\begin{aligned}
\T_s(K)\setminus u
=\T_s(\Cl(H-u)),~
\lk(\T_s(K),u)
=\T_{s-1}(\Cl(H-u))\cup\Sk_s(\rho).
\end{aligned}
\]
\end{lemma}

\begin{proof}
This follows from Lemma~\ref{lem:Tdl},
since $\lk(\Cl(H),u)=2^\rho$.
\end{proof}

The link in Lemma~\ref{lem:simplicial}
need not be a tolerance complex of a
clique complex. We therefore introduce
a larger class closed under the required
deletion and link operations.

\begin{definition}
A \emph{configuration}
$\mathfrak{C}=(H,s;\cS)$
 consists of a chordal
graph $H$ on $W$, an integer $s\ge-1$,
and a system $\cS$ whose terms $(\tau,r)$
satisfy: $\tau$ is a clique of $H$
and $r\ge s+1$. We require $\cS$ to be
nonempty if $s=-1$.

Its complex is
$Y(\mathfrak{C})=\T_s(\Cl(H))\cup Y(\cS)$.
The \emph{cost} of a chain is
$\cost(\cK)=\topp(\cK)-s\ge1$,
and the \emph{potential} of $\mathfrak{C}$ is
\[
\Phi(\mathfrak{C})=2s+1+
\sum_{\cK\in\cS}\cost(\cK).
\]
\end{definition}

Thus $Y(H,t;\emptyset)=\T_t(\Cl(H))$
and $\Phi(H,t;\emptyset)=2t+1$.
If $s=-1$, then $Y(\mathfrak{C})=Y(\cS)$
and $\Phi(\mathfrak{C})=\Psi(\cS)$.
Clean-up preserves the configuration
conditions and $Y(\mathfrak{C})$, does not
increase $\Phi$, and cannot remove
every term of a nonempty system.

A configuration is \emph{trivial} if
either $s\ge0$ and $H$ is complete,
or some term satisfies
$|W\setminus\tau|\le r$.
In either case, $Y(\mathfrak{C})=2^W$.

\begin{definition}
Let $\mathfrak{C}=(H,s;\cS)$ with $s\ge0$,
let $u$ be a simplicial vertex of $H$,
and put $\rho=N_H(u)$.
Let $\cK_u$ be a selected chain of
$\cS$, or let $\cK_u=\emptyset$,
meaning that no chain is selected.
Set $\cost(\emptyset)=0$.

The pair $(u,\cK_u)$ is
\emph{admissible} if
\begin{enumerate}[label=\textup{(A\arabic*)}]
\item \label{A1}every term of $\cS$ whose center
contains $u$ belongs to the selected
chain $\cK_u$;
\item \label{A2} every center of $\cK_u$ not
containing $u$ is contained in $\rho$.
\end{enumerate}
\end{definition}

For an admissible pair, define two
configurations on $W\setminus\{u\}$.

The \emph{deletion configuration} is
$\mathfrak{C}\setminus u=(H-u,s;\cS\setminus u)$,
where each center $\tau$ is replaced by
$\tau\setminus\{u\}$.

The \emph{link configuration} is
$\mathfrak{C}/u=(H-u,s-1;\cS/u)$.
For each unselected chain, replace
every term $(\tau,r)$ by $(\tau,r-1)$.
Replace the selected chain $\cK_u$
by a chain $\cK'$ obtained from the list
\[
\bigl((\tau,r-1):u\notin\tau\bigr),~
\bigl((\tau\setminus\{u\},r):u\in\tau\bigr),
~ ((\rho,s)),
\]
where the terms from $\cK_u$ retain
their original relative order.
At the junction of the first two
segments, perform the clean-up from
the proof of Theorem~\ref{thm:skeleta}: delete the terms $(\tau,r-1)$ with $u\notin\tau$ whose radius $r$ equals the radius of the first term of $\cK_u$ whose center contains $u$.
If $\cK_u=\emptyset$, then
$\cK'=((\rho,s))$.

\begin{lemma}\label{lem:step}
Let $(u,\cK_u)$ be admissible for
$\mathfrak{C}$ with $s\ge0$. Then
$u\in V(Y(\mathfrak{C}))$, both $\mathfrak{C}\setminus u$
and $\mathfrak{C}/u$ are configurations, and
\[
\begin{aligned}
Y(\mathfrak{C})\setminus u=Y(\mathfrak{C}\setminus u),~
\lk(Y(\mathfrak{C}),u)=Y(\mathfrak{C}/u).
\end{aligned}
\]
Moreover,
$\Phi(\mathfrak{C}\setminus u)=\Phi(\mathfrak{C})$
and $\Phi(\mathfrak{C}/u)\le\Phi(\mathfrak{C})-1$.
\end{lemma}

\begin{proof}
Since $s\ge0$, we have
$\{u\}\in\Cl(H)\subseteq Y(\mathfrak{C})$.
By Lemmas~\ref{lem:simplicial}
and~\ref{lem:Sdl}(1),
$Y(\mathfrak{C})\setminus u=Y(\mathfrak{C}\setminus u)$.
The graph $H-u$ is chordal, the
centers remain nested cliques, and
the radii are unchanged. Thus
$\mathfrak{C}\setminus u$ is a configuration
and $\Phi(\mathfrak{C}\setminus u)=\Phi(\mathfrak{C})$. 
By Lemmas~\ref{lem:simplicial}
and~\ref{lem:Sdl}(2), the complex
$\lk(Y(\mathfrak{C}),u)$ is the union of
$\T_{s-1}(\Cl(H-u))$, $\Sk_s(\rho)$,
and the balls
$\Sk_r(\tau\setminus\{u\})$
for terms with $u\in\tau$ and
$\Sk_{r-1}(\tau)$ for terms with
$u\notin\tau$.

By \ref{A1}, every unselected chain
contributes only terms of the latter
type, and its radii decrease by one.

In the selected chain, the centers
containing $u$ form a final segment.
After removing $u$, these centers
are contained in $\rho$, since they
are cliques containing $u$.
The preceding centers are also
contained in $\rho$ by $\ref{A2}$.
Thus the transformed centers remain
nested, and $(\rho,s)$ can be
appended. All radii are at least
$s$, and the only possible failure
of monotonicity occurs at the
junction of the first two segments.
The clean-up from
Theorem~\ref{thm:skeleta} removes
this failure.
All resulting centers are cliques
of $H-u$, and all radii are at
least $(s-1)+1$. Since $\cK'$
contains $(\rho,s)$, the system
$\cS/u$ is nonempty. Hence
$\mathfrak{C}/u$ is a configuration and
$\lk(Y(\mathfrak{C}),u)=Y(\mathfrak{C}/u)$.

For each unselected chain, both
the top and $s$ decrease by one,
so its cost is unchanged.
If $\cK_u\ne\emptyset$, the top
of $\cK'$ is at most
$\topp(\cK_u)$. If
$\cK_u=\emptyset$, its top is $s$.
In either case,
$\cost(\cK')\le\cost(\cK_u)+1$.
Consequently,
$\Phi(\mathfrak{C}/u)
\le2(s-1)+1+
\sum_{\cK\in\cS}\cost(\cK)+1
=\Phi(\mathfrak{C})-1.
$
\end{proof}

\begin{theorem}\label{thm:induction}
Let $\fF$ be a class of configurations such that, for
every nontrivial $\mathfrak{C}\in\fF$ with $s\ge0$, there is
an admissible pair $(u,\cK_u)$ for which some clean-ups
of $\mathfrak{C}\setminus u$ and $\mathfrak{C}/u$ belong to $\fF$.
Then $C(Y(\mathfrak{C}))\le\Phi(\mathfrak{C})$ for every
$\mathfrak{C}\in\fF$.
\end{theorem}

\begin{proof}
We induct on $|W|$. If $s=-1$, then
$\cS\ne\emptyset$, and
Theorem~\ref{thm:skeleta} gives
$C(Y(\mathfrak{C}))\le\Psi(\cS)=\Phi(\mathfrak{C})$.
If $s\ge0$ and $\mathfrak{C}$ is trivial, then
$C(Y(\mathfrak{C}))=0$.
Otherwise, choose an admissible pair
$(u,\cK_u)$ as in the hypothesis.
The cleaned-up link configuration has
a nonempty system, since it contains
$(\rho,s)$ before clean-up and clean-up
cannot remove every term.
The deletion configuration has $s\ge0$.
Thus, by induction and
Lemma~\ref{lem:step},
$C(Y(\mathfrak{C})\setminus u)\le\Phi(\mathfrak{C}),$
$C(\lk(Y(\mathfrak{C}),u))\le\Phi(\mathfrak{C})-1.$
Lemma~\ref{lem:tancer} now gives
$C(Y(\mathfrak{C}))\le\Phi(\mathfrak{C})$.
\end{proof}

\begin{corollary}\label{cor:induction}
If $H$ is chordal and $(H,t;\emptyset)$
belongs to a class $\fF$ as in
Theorem~\ref{thm:induction}, then
$\T_t(\Cl(H))$ is $(2t+1)$-collapsible.
\end{corollary}

\begin{figure}[ht]
\centering
\begin{tikzpicture}[scale=0.9,
  node/.style={circle,fill=black,inner sep=1.2pt},
  sub/.style={line width=1.6pt},
  ring/.style={circle,draw,line width=1.1pt,inner sep=3.2pt},
  lab/.style={font=\small}]
  % Host spider
  \coordinate (q) at (0,0);
  \coordinate (a1) at (0,1.2);
  \coordinate (a2) at (0,2.4);
  \coordinate (b1) at (-1.05,-0.6);
  \coordinate (b2) at (-2.1,-1.2);
  \coordinate (c1) at (1.05,-0.6);
  \coordinate (c2) at (2.1,-1.2);
  \draw (a2)--(q)--(b2) (q)--(c2);
  \foreach \p in {q,a1,a2,b1,b2,c1,c2}
    \node[node] at (\p) {};
  \node[lab,right] at (0.12,0.1) {$q$};

  % Subtrees
  \node[ring,label={[lab]above:$a$}] at (a2) {};
  \draw[sub] (0.22,1.15)--(0.22,2.35);
  \node[lab,right] at (0.24,1.75) {$b$};
  \draw[sub] (-0.22,1.2)--(-0.22,0.12)--(-1.15,-0.42);
  \node[lab,left] at (-0.24,0.65) {$c$};
  \draw[sub] (-0.93,-0.82)--(-1.98,-1.42);
  \node[lab,below] at (-1.35,-1.1) {$e$};
  \node[ring,label={[lab]left:$f$}] at (b2) {};
  \draw[sub] (-0.02,-0.28)--(1.98,-1.44);
  \node[lab,below] at (1.25,-0.95) {$g$};
  \node[ring,label={[lab]above right:$h$}] at (c1) {};

  % Intersection graph
  \begin{scope}[xshift=6.2cm,yshift=0.2cm]
    \coordinate (C) at (0,0);
    \coordinate (B) at (0,1.05);
    \coordinate (A) at (0,2.1);
    \coordinate (E) at (-0.95,-0.55);
    \coordinate (F) at (-1.9,-1.1);
    \coordinate (G) at (0.95,-0.55);
    \coordinate (H) at (1.9,-1.1);
    \draw (A)--(B)--(C)--(E)--(F)
          (C)--(G)--(H);
    \foreach \p in {A,B,C,E,F,G,H}
      \node[node] at (\p) {};
    \node[lab,right] at (A) {$a$};
    \node[lab,right] at (B) {$b$};
    \node[lab,right] at (0.1,0.12) {$c$};
    \node[lab,above left] at (E) {$e$};
    \node[lab,above left] at (F) {$f$};
    \node[lab,above right] at (G) {$g$};
    \node[lab,above right] at (H) {$h$};
  \end{scope}
\end{tikzpicture}
\caption{Left: a spider $\Gamma$ with center $q$
and seven subtrees, indicated by thick segments
and circles beside its vertices. Right: their
intersection graph $H$, which is neither an
interval graph nor a split graph. Here $\Gamma_a$ is the top vertex, $\Gamma_b$ the upper edge of the vertical leg, $\Gamma_c$ the path from the middle vertex of the vertical leg through $q$ to the middle vertex of the lower-left leg, $\Gamma_e$ the outer edge of the lower-left leg, $\Gamma_f$ its endpoint, $\Gamma_g$ the lower-right leg together with $q$, and $\Gamma_h$ the middle vertex of the lower-right leg.}
\label{fig:spider}
\end{figure}
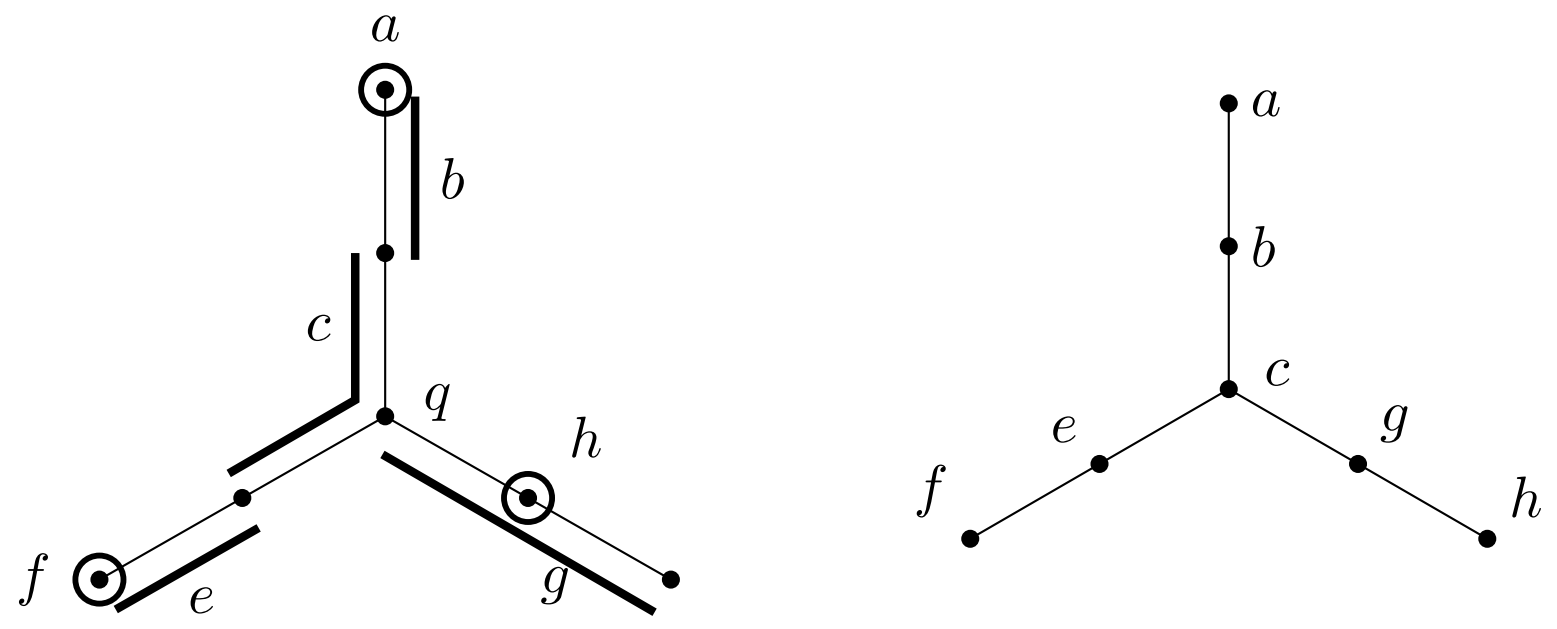

We first apply Theorem~\ref{thm:induction}
to chordal graphs defined by a host tree.
A \emph{spider} is a finite tree with at most one
vertex of degree at least $3$.
A graph $H$ on $W$ is a \emph{subtree graph of a spider}
if there are a spider $\Gamma$ and nonempty
subtrees $(\Gamma_w)_{w\in W}$ such that
distinct vertices $u,v$ are adjacent if
and only if $\Gamma_u\cap\Gamma_v\ne\emptyset$.
These graphs are chordal \cite{Gav74} and are closed
under taking induced subgraphs.
Interval graphs arise when $\Gamma$ is a
path; Figure~\ref{fig:spider} gives an
example that is not an interval graph.

Split graphs are also subtree graphs of spiders.
Indeed, if $W=C\sqcup I$, where $C$ is
a clique and $I$ is independent, take
$\Gamma$ to be a star with center $q$
and leaves $\ell_i$, $i\in I$.
Set $\Gamma_i=\{\ell_i\}$ for $i\in I$,
and for $c\in C$ let $\Gamma_c$ be the
subtree induced by
$\{q\}\cup\{\ell_i:ci\in E(H)\}$.

\begin{theorem}\label{thm:spider}
Let $H$ be a subtree graph of a spider and let $t\ge0$.
Then $\T_t(\Cl(H))$ is $(2t+1)$-collapsible.
Consequently, if $K$ is a complex on a finite
set $W$ with $K=\Cl(H)$ on its vertex set,
then
$C(\T_t(K))\le2t+1=\eta(2,t+1)-1$.
In particular, the case $d=1$ of
Conjecture~\ref{conj:coll} holds for
nerves of finite families of intervals
in $\mathbb R$.
\end{theorem}

\begin{proof}
Fix a vertex $q$ of $\Gamma$ of degree
at least $3$, or any vertex if none exists; in the latter case $\Gamma$ is a path and there are at most two legs.
The components of $\Gamma-q$ are paths,
called \emph{legs}, each attached to $q$
at one end. Let $\delta(p)$ denote the
distance from $q$ to a vertex $p$.
If $q\notin\Gamma_w$, then $\Gamma_w$
lies in a unique leg, and
$\operatorname{apex}(w)$ denotes its
vertex of minimum distance from $q$.
For $U\subseteq W$ and a vertex $p$
of $\Gamma$, put
$U_p=\{w\in U:p\in\Gamma_w\}$,
which is a clique of $H$.
We use the following observation,
which follows from the connectedness
of the subtrees and the fact that
each leg is a path.

\smallskip
\noindent
$(\star)$
Let $L$ be a leg, let $x,p\in L$ with
$\delta(p)\ge\delta(x)$, and suppose
$p\in\Gamma_w$. If $q\in\Gamma_w$,
or if $\Gamma_w\subseteq L$ and
$\delta(\operatorname{apex}(w))\le\delta(x)$,
then $x\in\Gamma_w$.
\smallskip

Let $\fF_{\mathrm{sp}}$ be the class of configurations $(H[U],s;\cS)$, where
$U\subseteq W$, whose chains are assigned to pairwise distinct legs and
satisfy the following condition: every center $\tau$ of a chain assigned to
a leg $L$ is contained in $U_p$ for some $p\in L$ such that
$\delta(p)\ge\delta(\operatorname{apex}(w))$ for every $w\in U$ with
$\Gamma_w\subseteq L$.
Clearly $(H,t;\emptyset)\in\fF_{\mathrm{sp}}$.
Let $\mathfrak{C}\in\fF_{\mathrm{sp}}$ be nontrivial with $s\ge0$. Since $H[U]$ is not
complete, $U\ne U_q$, so some $w\in U$ satisfies $q\notin\Gamma_w$; its
subtree $\Gamma_w$ then lies in a leg. Among all such $w$, over all legs,
choose $u$ with $\delta(\operatorname{apex}(u))$ maximum. Let $L$ be the leg
containing $\Gamma_u$, and put $x=\operatorname{apex}(u)$. Let $\cK_u$ be the
chain assigned to $L$, or $\emptyset$ if no such chain exists.

We claim that
$N_{H[U]}(u)=U_x\setminus\{u\}$.
Indeed, if $w$ is adjacent to $u$, then
$\Gamma_w$ contains a vertex $p$ of
$\Gamma_u$, with $\delta(p)\ge\delta(x)$.
By $(\star)$ and the maximality of
$\delta(x)$, we have $x\in\Gamma_w$.
The reverse inclusion is immediate.
Thus $u$ is simplicial, with
$\rho=U_x\setminus\{u\}$.
The pair $(u,\cK_u)$ is admissible.
For \ref{A1}, every center of a chain assigned
to another leg is contained in some $U_p$
with $p\notin L$, whereas
$\Gamma_u\subseteq L$.
For $\ref{A2}$, let $\tau\subseteq U_p$
be a center of $\cK_u$ not containing
$u$, and let $w\in\tau$.
The defining condition gives
$\delta(p)\ge\delta(x)$.
By $(\star)$ and the choice of $u$,
$x\in\Gamma_w$, so $w\in\rho$.
Hence $\tau\subseteq\rho$.

Clean-ups of $\mathfrak{C}\setminus u$ and
$\mathfrak{C}/u$ belong to $\fF_{\mathrm{sp}}$.
Indeed, existing centers only lose
vertices, and the defining conditions
weaken as $U$ shrinks. The new term
$(\rho,s)$ joins the chain assigned
to $L$, or forms that chain, and
$\rho\subseteq U_x$ with $\delta(x)$
at least the apex distances of all
remaining subtrees contained in $L$.
Clean-up only removes terms.
Theorem~\ref{thm:induction} and
Corollary~\ref{cor:induction} now give
$C(\T_t(\Cl(H)))\le2t+1$.
The bound holds for every $0\le s\le t$,
so the second assertion follows from
Lemma~\ref{lem:ghost}. Finally, by
Helly's theorem on the line, the nerve
of a finite family of intervals is,
on its vertex set, the clique complex
of an interval graph.
\end{proof}

The bound $2t+1$ is sharp by
Example~\ref{ex:sharp} with $m=2$.
The argument imposes no bound on the
number of chains: the selected vertex
never belongs to centers assigned to
two different legs.

We next show that the induction principle
applies to all chordal graphs when $t\le3$.
The key observation is that condition \ref{A2}
is automatic for a chain whose last center
contains the selected vertex.

\begin{lemma}\label{lem:topadmissible}
Let $\mathfrak{C}=(H,s;\cS)$ be a configuration
with $s\ge0$, and let $u$ be a simplicial
vertex of $H$ belonging to the last center
of at most one chain of $\cS$.
Let $\cK_u$ be that chain, if it exists,
and $\cK_u=\emptyset$ otherwise.
Then $(u,\cK_u)$ is admissible.
\end{lemma}

\begin{proof}
Since the centers of each chain are nested,
a chain has a center containing $u$ only
if its last center contains $u$.
This gives \ref{A1}. For \ref{A2}, let $A$ be
the last center of $\cK_u$.
Since $A$ is a clique containing $u$,
$A\setminus\{u\}\subseteq N_H(u)$.
Every center of $\cK_u$ not containing
$u$ is contained in $A\setminus\{u\}$.
\end{proof}

\begin{theorem}\label{thm:smallt}
Let $H$ be a chordal graph and
$0\le t\le3$. Then $\T_t(\Cl(H))$
is $(2t+1)$-collapsible.
Consequently, for every $1$-collapsible
complex $K$ and every $t\ge0$ with
$t\le3$ or $t\ge\dim K$,
\[
C(\T_t(K))\le2t+1=\eta(2,t+1)-1.
\]
Thus the case $d=1$ of
Conjecture~\ref{conj:coll} holds for
$t\le3$, and for all $t$ when
$\dim K\le4$.
\end{theorem}

\begin{proof}
Let $\fF_3$ be the class of configurations
$(H,s;\cS)$ satisfying $s+|\cS|\le3$,
where $|\cS|$ denotes the number of chains.
It contains $(H,t;\emptyset)$ for $t\le3$.
The class is closed under the operations
of Lemma~\ref{lem:step} and clean-up:
deletion preserves $s$ and the number
of chains, while the link decreases $s$
by one and adds at most one chain.
Clean-up cannot increase the number
of chains.
Let $\mathfrak{C}=(H,s;\cS)\in\fF_3$ be
nontrivial with $s\ge0$. Then
$|\cS|\le3$ and $H$ is not complete.
By Dirac's theorem, $H$ has two
nonadjacent simplicial vertices $x,y$.
Since the last center of each chain
is a clique, it cannot contain both
$x$ and $y$. Thus the number of chains whose last center contains $x$ plus the number whose last center contains $y$ is at most $|\cS|\le3$, so one of $x$, $y$ lies in at most one last center. By
Lemma~\ref{lem:topadmissible}, it gives
an admissible pair. Theorem~\ref{thm:induction}
and Corollary~\ref{cor:induction} prove
the first assertion.

For the second, let $K$ be a
$1$-collapsible complex on $W$.
The void case is immediate, so assume
$K$ is nonvoid. If $t\ge\dim K$, then
$C(\T_t(K))\le\dim K+1+t\le2t+1$
by Remark~\ref{rem:dimbound}.
If $t\le3$, Lemma~\ref{lem:chordal}
gives $K=\Cl(H)$ on its vertex set,
with $H$ chordal. The first assertion
and Lemma~\ref{lem:ghost} apply.
Finally, if $\dim K\le4$, every
$t\ge0$ satisfies $t\le3$ or
$t\ge\dim K$.
\end{proof}

For $t=0$, the invariants $C(\Cl(H))$, $L_\kk(\Cl(H))$, and $h(\Cl(H))$
coincide by Lemma~\ref{lem:chordal}. The same holds for $\T_t(\Cl(H))$ when
$\Cl(H)$ is shifted, equivalently when $H$ is a threshold graph, by
Theorem~\ref{thm:stable} applied after reversing the vertex order.
Since $I_{\Cl(H)}=I(\overline H)$, one might therefore expect that the
regularity of $\sym{I(G)}{k}$ always equals the largest degree of a minimal
generator when $\overline G$ is chordal, that is, when $G$ is co-chordal.
The following construction shows that
this fails, with arbitrarily large gaps.

\begin{theorem}\label{thm:splitgap}
Let $\Delta$ be a nonvoid complex on a
finite set $Q$ with $\Delta\ne2^Q$, and let
$A_0,\dots,A_t$ be its facets, where $t\ge0$.
Let $H$ be the split graph consisting
of a clique on $Q$ and, for each $i$,
a new vertex $z_i$ adjacent exactly to
the vertices of $A_i$.
Put $m=x_{z_0}\cdots x_{z_t}$.
Then
$I_{\T_t(\Cl(H))}=m\,I_\Delta$.
Consequently, for every field $\kk$,
\[
h(\T_t(\Cl(H)))=t+1+h(\Delta),
~
L_\kk(\T_t(\Cl(H)))=t+1+L_\kk(\Delta).
\]
\end{theorem}
 
\begin{proof}
Let $Z=Y\sqcup B$, where $Y\subseteq Q$
and $B\subseteq\{z_0,\dots,z_t\}$.
Since the faces of $\Cl(H)$ are the
cliques of $H$,
Proposition~\ref{prop:gens} gives
\[
Z\notin\T_t(\Cl(H))
\iff |Z|-\omega(H[Z])\ge t+1,
\]
where $\omega$ denotes the clique number,
with $\omega(H[\emptyset])=0$.
 
Since $\omega(H[Z])\ge|Y|$ and
$|B|\le t+1$, this requires
$|B|=t+1$ and $\omega(H[Z])=|Y|$.
When $B$ contains every $z_i$,
the latter condition holds precisely
when $Y$ is contained in no $A_i$.
Indeed, a largest clique containing $z_i$
has exactly $1+|Y\cap A_i|$ vertices,
whereas $Y$ itself is a clique.
Finally, $Y$ is contained in no $A_i$
if and only if $Y\notin\Delta$.
Thus $I_{\T_t(\Cl(H))}=m\,I_\Delta$. 
The minimal generators of $I_\Delta$ have
largest degree $h(\Delta)+1$, and
$\deg m=t+1$. This gives the formula for
the Helly dimension. Multiplication by $m$
shifts the regularity of the quotient by
$t+1$, so the formula for the Leray number
follows from Proposition~\ref{prop:sym}
and Hochster's formula.
\end{proof}

\begin{theorem}\label{thm:notHelly}
Let $r,c\ge2$, let $P_1,\dots,P_r$
be pairwise disjoint sets of size $c$,
and put $Q=P_1\sqcup\dots\sqcup P_r$.
Let $\mathcal A$ be the family of
$c^r$ subsets of $Q$ obtained by
deleting one element from each $P_i$.
Let $H$ be the chordal graph consisting
of a clique on $Q$ and, for each
$A\in\mathcal A$, a new vertex $z_A$
adjacent exactly to the vertices of $A$.
Put $t=c^r-1$ and
$m=\prod_{A\in\mathcal A}x_{z_A}$.
Then
$I_{\T_t(\Cl(H))}
=m(x^{P_1},\dots,x^{P_r}).
$
Consequently, for every field $\kk$,
\[
h(\T_t(\Cl(H)))=c^r+c-1,
~
L_\kk(\T_t(\Cl(H)))=c^r+r(c-1).
\]
Moreover, for every integer $t\ge3$,
there is a chordal graph $H$ such that
$L_\kk(\T_t(\Cl(H)))>h(\T_t(\Cl(H)))$.
\end{theorem}
 
\begin{proof}
Apply Theorem~\ref{thm:splitgap} to
$\Delta=\partial2^{P_1}*\dots*\partial2^{P_r}$.
A subset of $Q$ is a nonface of $\Delta$ if and
only if it contains some $P_i$. Hence the facets
of $\Delta$ are the $c^r$ sets in $\mathcal A$,
and $I_\Delta=(x^{P_1},\dots,x^{P_r})$.
Thus $h(\Delta)=c-1$. Since the monomials
$x^{P_i}$ form a regular sequence of degree $c$,
$L_\kk(\Delta)=\reg S/(x^{P_1},\dots,x^{P_r})
=r(c-1)$. This gives the stated formulas.
In particular, the difference is $(r-1)(c-1)$.

For the last assertion, start with
$r=c=2$, so $t=3$. For this graph,
$|\sigma|-\omega(H[\sigma])\le t+1$
for all $\sigma\subseteq V(H)$:
if $\sigma=Y\sqcup B$ as in the proof of
Theorem~\ref{thm:splitgap}, then
$\omega(H[\sigma])\ge|Y|$ and $|B|\le t+1$.
More generally,
suppose $H$ satisfies this inequality,
and let $H'=H\sqcup\{z\}$, where
$z$ is isolated. Then
$I_{\T_{t+1}(\Cl(H'))}
=x_z I_{\T_t(\Cl(H))}.
$ Indeed, by Proposition~\ref{prop:gens} a set $\sigma'\subseteq V(H')$ is a nonface of $\T_{t+1}(\Cl(H'))$ if and only if $|\sigma'|-\omega(H'[\sigma'])\ge t+2$. If $z\notin\sigma'$, this fails by the inequality. If $\sigma'=\sigma\cup\{z\}$, then it fails for $\sigma=\emptyset$, and for $\sigma\ne\emptyset$ it says $|\sigma|-\omega(H[\sigma])\ge t+1$, that is, $\sigma\notin\T_t(\Cl(H))$.
Thus adjoining an isolated vertex
and increasing the tolerance by one
increases both the Helly dimension
and the Leray number by one.
The same inequality on clique
numbers holds for $H'$ with
$t$ replaced by $t+1$, so the
construction can be iterated.
Starting with $r=c=2$ gives
$L_\kk>h$ for every $t\ge3$.
\end{proof}
 
The smallest instance of this construction, with $r=c=2$, has eight
vertices and $t=3$. Here $H$ is $K_4$ with one new vertex attached to each
edge of a $4$-cycle, and
$h(\T_3(\Cl(H)))=5$ while $L_\kk(\T_3(\Cl(H)))=6$.

A complex $\Delta$ with at most three facets satisfies
$L_\kk(\Delta)=h(\Delta)$, so Theorem~\ref{thm:splitgap} gives a gap only
when $t\ge3$. Indeed, by Corollary~\ref{cor:facets} and
Lemma~\ref{lem:hLC}, $L_\kk(\Delta)\le C(\Delta)\le2$.
If $L_\kk(\Delta)=0$, then $\Delta$ is a simplex on its vertex set, and
$h(\Delta)=0$. If $L_\kk(\Delta)=1$, then $\Delta$ is not a simplex on its
vertex set, so $h(\Delta)\ge1$, and $h(\Delta)\le L_\kk(\Delta)=1$ by
Lemma~\ref{lem:hLC}.
Suppose $L_\kk(\Delta)=2$, and choose $U$ with
$\Ht_1(\Delta[U];\kk)\ne0$. The facets of $\Delta[U]$ are among the sets
$F\cap U$, where $F$ ranges over the facets of $\Delta$, so there are at
most three of them. By the nerve theorem \cite[Theorem~10.6]{Bjo95}, $\Delta[U]$ is homotopy equivalent
to the nerve of its facets, which therefore has nonzero first homology.
Hence $\Delta$ has exactly three facets $F_1,F_2,F_3$, the sets
$F_i\cap U$ are the facets of $\Delta[U]$, they intersect pairwise, and
$F_1\cap F_2\cap F_3\cap U=\emptyset$. Choose
$p_{ij}\in F_i\cap F_j\cap U$ for $1\le i<j\le3$. These three points are
distinct. The set $\{p_{12},p_{13},p_{23}\}$ is not contained in any $F_k$,
since otherwise one of the points would lie in $F_1\cap F_2\cap F_3\cap U$.
Each pair, however, lies in a common facet. Thus
$\{p_{12},p_{13},p_{23}\}$ is a missing face of $\Delta$, and
$h(\Delta)=2$.

For $t=1$, equality holds for every chordal graph.
 
\begin{proposition}\label{prop:t1}
For every chordal graph $H$ and every field $\kk$,
\[
C(\T_1(\Cl(H)))=L_\kk(\T_1(\Cl(H)))
=h(\T_1(\Cl(H))).
\]
\end{proposition}
 
\begin{proof}
Suppose first that $H$ has an induced subgraph
on a set $\sigma$ of four vertices that is
$2K_2$ or $P_4$. Then $\omega(H[\sigma])=2$,
so $\sigma\notin\T_1(\Cl(H))$. Every three-element
subset of $\sigma$ contains an edge, so it lies
in $\T_1(\Cl(H))$. Thus $\sigma$ is a missing
face, and Lemma~\ref{lem:hLC} and
Theorem~\ref{thm:smallt} give
$3\le h(\T_1(\Cl(H)))\le L_\kk(\T_1(\Cl(H)))
\le C(\T_1(\Cl(H)))\le3.
$
Otherwise, $H$ has no induced $2K_2$ or $P_4$.
Since $H$ is chordal, it also has no induced $C_4$, so $H$ is a threshold
graph \cite{CH77}. Hence $V(H)$ can be ordered so that $u<v$ and
$v$ adjacent to $w\ne u$ imply $u$ adjacent to $w$; order the vertices by nonincreasing degree, since threshold graphs have nested neighborhoods \cite{MP95}. If $f$ is a clique,
$x\in f$, $y<x$ and $y\notin f$, then $y$ is adjacent to every vertex of
$f\setminus\{x\}$, so $(f\setminus\{x\})\cup\{y\}$ is a clique. Thus
$\Cl(H)$ is shifted, and Theorem~\ref{thm:stable} applies after reversing
the vertex order.
\end{proof}
 
We do not know whether
$L_\kk(\T_2(\Cl(H)))=h(\T_2(\Cl(H)))$
for all chordal graphs $H$; a computer check found equality for every chordal graph with at most seven vertices, over $\mathbb F_2$ and over a large prime field. 
By contrast, the invariants $C$ and $L_\kk$ of
$\T_t(\Cl(H))$ differ in general, because
the Leray number can depend on the field.
 
\begin{corollary}\label{cor:field}
For every field $\kk$, there are a chordal graph
$H$ and an integer $t\ge3$ such that
$C(\T_t(\Cl(H)))>L_\kk(\T_t(\Cl(H)))$.
\end{corollary}
 
\begin{proof}
Choose a prime $p$ different from the
characteristic of $\kk$. Let $\Delta$ be a finite
two-dimensional simplicial complex with
$\widetilde H_1(\Delta;\mathbb Z)\cong\mathbb Z/p$ and
$\widetilde H_2(\Delta;\mathbb Z)=0$, for instance a
triangulation of the mapping cone of a degree-$p$
map of the circle; for $p=2$ one can take the six-vertex real projective plane of Example~\ref{ex:RP2}. Let $t+1$ be its number of facets.
Then $\widetilde H_2(\Delta;\kk)=0$. Since $\dim\Delta=2$,
no induced subcomplex has a nonzero $2$-cycle
over $\kk$, so $L_\kk(\Delta)\le2$.
On the other hand,
$\widetilde H_2(\Delta;\mathbb F_p)\cong\mathbb F_p$,
so $L_{\mathbb F_p}(\Delta)=3$.
A nonzero $2$-cycle has at least four triangles
in its support, since each edge of the support lies in at least two triangles of the support, so $t\ge3$.
Let $H$ be the graph of
Theorem~\ref{thm:splitgap} for $\Delta$, and put
$T=\T_t(\Cl(H))$. Since $C(T)$ does not depend
on the field, Lemma~\ref{lem:hLC} gives
$C(T)\ge L_{\mathbb F_p}(T)=t+4>t+3\ge L_\kk(T)$.
\end{proof}
 
By Lemma~\ref{lem:chordal} and
Proposition~\ref{prop:t1}, no such example
exists for $t\le1$. We do not know the
smallest $t$ for which one exists.
 
\begin{example}\label{ex:RP2}
Let $\Delta$ be the six-vertex triangulation of the real projective plane,
with facets $\{1,2,4\}$, $\{1,2,6\}$, $\{1,3,4\}$, $\{1,3,5\}$, $\{1,5,6\}$, $\{2,3,5\}$, $\{2,3,6\}$, $\{2,4,5\}$, $\{3,4,6\}$,
and $\{4,5,6\}$. Then $t=9$, and the graph $H$ of Theorem~\ref{thm:splitgap} has
$16$ vertices. Every pair of vertices is an edge of $\Delta$, so the missing
faces of $\Delta$ are the ten triangles that are not facets, and
$h(\Delta)=2$. Hence $L_\kk(\Delta)\ge2$ by Lemma~\ref{lem:hLC}, while
$L_\kk(\Delta)\le3$ since $\dim\Delta=2$. A nonzero $2$-cycle of an induced
subcomplex is a nonzero $2$-cycle of $\Delta$; as
$\Ht_2(\Delta;\kk)=0$ for $\operatorname{char}\kk\ne2$ and
$\Ht_2(\Delta;\mathbb F_2)\cong\mathbb F_2$, we get $L_\kk(\Delta)=2$ in the
first case and $L_{\mathbb F_2}(\Delta)=3$. By Theorem~\ref{thm:splitgap},
\[
L_{\mathbb F_2}(\T_9(\Cl(H)))=13,~
L_\kk(\T_9(\Cl(H)))=h(\T_9(\Cl(H)))=12
~\text{whenever }\operatorname{char}\kk\ne2.
\]
\end{example}

\section{The case \texorpdfstring{$t=1$}{t=1}}
\label{sec:tone}

In this section, we reduce the conjectured Leray bound for the $1$-tolerance complexes of
$d$-collapsible complexes to a single inequality involving intersections of links
(Question~\ref{q:links}). The method refines the argument of \cite[Theorem~1.5]{KL23},
which bounds the
homology of a union of links using the Kalai--Meshulam union theorem.

\begin{lemma}\label{lem:Rhelly}
Let $\Lambda$ be an $e$-Leray complex on
$T\sqcup U$, where $|T|=r\ge1$ and
$T\notin\Lambda$, and put
$R=\bigcap_{v\in T}
\lk(\Lambda,T\setminus\{v\})[U]$.
If $R$ is nonvoid, then every missing face
of $R$ has cardinality at most $e+2-r$.
\end{lemma}

\begin{proof}
Let $\xi$ be a missing face of $R$.
Since $R$ is nonvoid, $\xi\ne\emptyset$.
Put
$c=\{A\subsetneq T:A\cup\xi\in\Lambda\}$.
Since $\xi\notin R$, the complex $c$ is
a proper, possibly void, subcomplex of $\partial2^T$; for $r=1$, $\partial2^T=\{\emptyset\}$ and $c$ is void.
Moreover, $R$ is nonvoid and every proper
subset of $\xi$ belongs to $R$.
Together with $T\notin\Lambda$, this gives
\[
\Lambda[T\cup\xi]
=(\partial2^T*\partial2^\xi)
\cup(c*2^\xi).
\]

The first complex is a homology sphere of
dimension $m=r+|\xi|-3$, and the second
is void or acyclic. Their intersection is
$c*\partial2^\xi$. By the K\"unneth
formula for joins,
\[
\Ht_m(c*\partial2^\xi;\kk)
\cong
\Ht_{r-2}(c;\kk)\otimes_\kk
\Ht_{|\xi|-2}(\partial2^\xi;\kk)=0,
\]
since $c$ is a proper subcomplex of the
$(r-2)$-sphere $\partial2^T$.
The Mayer--Vietoris sequence therefore gives $\widetilde H_m(\Lambda[T\cup\xi];\Bbbk)\ne0$.
If
$m\ge0$, then $m\le e-1$ because $\Lambda$ is $e$-Leray; if $m=-1$, then $r=|\xi|=1$ and
$m\le e-1$ holds since $e\ge0$. In either case $|\xi|\le e+2-r$.
\end{proof}

\begin{question}\label{q:links}
Let $\Lambda$ be an $e$-Leray complex on $T\sqcup U$, where $|T|=r\ge1$ and
$T\notin\Lambda$, and put
$R=\bigcap_{v\in T}\lk(\Lambda,T\setminus\{v\})[U]$.
Is $L_\kk(R)\le r(e+1-r)$ whenever $R$ is nonvoid?
\end{question}

If $R$ is nonvoid, then $T\setminus\{v\}\in\Lambda$ for every $v\in T$, so
$\Lambda[T]=\partial2^T$; since $\Lambda$ is $e$-Leray, this forces
$r\le e+1$. For $r=1$ the answer is positive, since then $R=\Lambda[U]$.
For $r=e+1$ it is also positive, since then $R$ is a simplex on its vertex
set by Lemma~\ref{lem:Rhelly}. More generally, the answer is positive
whenever the Stanley--Reisner ideals of the complexes
$\lk(\Lambda,T\setminus\{v\})[U]$, $v\in T$, are componentwise linear.
Indeed, by Lemma~\ref{lem:Rhelly}, $I_R$ is generated by the elements of
degree at most $q=e+2-r$ of these ideals; truncations of componentwise
linear ideals in degrees at most $q$ have regularity at most $q$
\cite{HH99}; and the Kalai--Meshulam bound \cite{KM06} then gives
$L_\kk(R)\le r(q-1)=r(e+1-r)$.
By Theorem~\ref{thm:t1conditional}, a positive answer to Question~\ref{q:links} would
establish the $t=1$ case of Conjecture~\ref{conj:leray} for $d$-collapsible
complexes. Since $1\le r\le e+1$, the cases above show that the answer is positive for $e\le1$. Hence Theorem~\ref{thm:t1conditional} holds unconditionally for $d\le2$, and for $d=2$ it recovers \cite[Theorem~1.6]{KL23}. The first open case is $e=r=2$: then $R$ is a flag complex by Lemma~\ref{lem:Rhelly}, and the question asks whether $L_\kk(R)\le2$. A positive answer would show that $\T_1(K)$ is $8$-Leray for every $3$-collapsible complex $K$.

\begin{theorem}\label{thm:t1conditional}
Let $d\ge1$, and suppose that Question~\ref{q:links} has a positive answer for $e=d-1$.
Then, for every $d$-collapsible complex $K$, the complex $\T_1(K)$ is
$(\eta(d+1,2)-1)$-Leray
over $\kk$; equivalently, $\reg(S/\sym{I_K}{2})\le\lfloor(d+3)^2/4\rfloor-1$.
\end{theorem}

\begin{proof}
Put $D=\lfloor(d+3)^2/4\rfloor-1$, and fix a sequence $K=K_0,\dots,K_n$ of elementary
$d$-collapses ending with a simplex or the void complex. Induced subcomplexes of
$d$-collapsible complexes are $d$-collapsible, so, by Lemma~\ref{lem:Tind}, it suffices to
prove that $\Ht_j(\T_1(K);\kk)=0$ for $j\ge D$.

We induct on $n$. If $K$ is a simplex or void, then
$L_\kk(\T_1(K))\le1\le D$. Otherwise, the first collapse is not at
$\emptyset$, since a collapse at $\emptyset$ removes every face and would
force $K$ to be a simplex. Let $\sigma\ne\emptyset$ be its free face, put
$a=|\sigma|\le d$, and write its unique maximal face as
$\tau=\sigma\sqcup U$. Let $K'=K_1$. By induction and the long exact
sequence of the pair, it suffices to show that
$H_j(\T_1(K),\T_1(K');\kk)=0$ $(j\ge D).$
If $U=\emptyset$, the two tolerance complexes differ only in faces of
cardinality at most $a+1$. Their relative homology therefore vanishes in
degrees greater than $a$. Since $a\le d<D$, the assertion follows.
Suppose $U\ne\emptyset$. By \cite[Proposition~4.3]{KL23}, applied to the elementary collapse $K\to K'$ at $\sigma$,
\[
H_j(\T_1(K),\T_1(K');\kk)
\;\cong
\bigoplus_{w\in W\setminus\tau}
\Ht_{j-a-1}\Bigl(
\bigcup_{v\in\sigma}\lk(K,\sigma\setminus\{v\})[U\cup\{w\}];\kk
\Bigr).
\]
Fix $w\in W\setminus\tau$. Since $U\cup(\sigma\setminus\{v\})\in K$ for
every $v\in\sigma$, the union inside the homology group is
$2^U\cup(w*B)$, where
\[
B=\bigcup_{v\in\sigma}B_v,
\qquad
B_v=\lk\bigl(K,(\sigma\setminus\{v\})\cup\{w\}\bigr)[U].
\]
If $B$ is void, the union is the simplex $2^U$, and the corresponding
summand vanishes. Otherwise, since $U\ne\emptyset$, the Mayer--Vietoris
sequence gives
$\Ht_{j-a-1}(2^U\cup(w*B);\kk)\cong\Ht_{j-a-2}(B;\kk)$.
We discard the void complexes $B_v$, and write $\sigma_0$ for the set of
$v\in\sigma$ with $B_v$ nonvoid; thus $B=\bigcup_{v\in\sigma_0}B_v$.

Let $\Lambda=\lk(K,w)[\tau]$. We claim that $\Lambda$ is $(d-1)$-Leray.
Indeed, for every nonempty $Z\subseteq\tau$,
$K[Z\cup\{w\}]=2^Z\cup(w*\Lambda[Z])$. Since $2^Z$ is a simplex, the
Mayer--Vietoris sequence gives
$\Ht_{\ell+1}(K[Z\cup\{w\}];\kk)\cong\Ht_\ell(\Lambda[Z];\kk)$ for
$\ell\ge0$. The same vanishing holds for $Z=\emptyset$ in these degrees.
Since $K$ is $d$-collapsible, it is $d$-Leray by Lemma~\ref{lem:hLC}, and
hence $\Lambda$ is $(d-1)$-Leray. Moreover, $\sigma\notin\Lambda$, since
$\sigma$ is contained only in the facet $\tau$.
Let $\emptyset\ne T\subseteq\sigma_0$, put $r=|T|$, and let
$R_T=\bigcap_{v\in T}B_v$. Each $B_v$ with $v\in T$ contains
$\emptyset$, so $R_T$ is nonvoid. Moreover,
$(\sigma\setminus T)\cup\{w\}\in K$, so $\sigma\setminus T\in\Lambda$, and
we may set $\Lambda_0=\lk(\Lambda,\sigma\setminus T)$, a complex on
$T\sqcup U$. Then $\Lambda_0$ is $(d-1)$-Leray, $T\notin\Lambda_0$, and
$R_T=\bigcap_{v\in T}\lk(\Lambda_0,T\setminus\{v\})[U].$

By the positive answer to Question~\ref{q:links}, applied to $\Lambda_0$ with $e=d-1$,
\begin{equation}\label{eq:RT} L_\kk(R_T)\le r(d-r). \end{equation}
% \begin{equation}\label{eq:RT}
% L_\kk(R_T)\le r(q-1)=r(d-r).
% \end{equation}

For complexes $X$ and $Y$ on the same ground set, the Mayer--Vietoris
sequence, applied to their induced subcomplexes, gives
$L_\kk(X\cup Y)\le\max\{L_\kk(X),\,L_\kk(Y),\,L_\kk(X\cap Y)+1\}$.
By induction on the number of complexes, it follows that, for complexes $X_1,\dots,X_n$ on one ground set, $L_\kk(\bigcup_iX_i)\le\max_{\emptyset\ne J}\bigl(L_\kk(\bigcap_{i\in J}X_i)+|J|-1\bigr)$ (cf.\ \cite{KM06}); in particular
\[
L_\kk(B)\le\max_{\emptyset\ne T\subseteq\sigma_0}
\bigl(L_\kk(R_T)+|T|-1\bigr).
\]
Indeed, writing $B=X\cup B_{v_0}$ with $X=\bigcup_{v\ne v_0}B_v$, the
intersection $X\cap B_{v_0}=\bigcup_{v\ne v_0}(B_v\cap B_{v_0})$ is again a
union of complexes whose intersections are the complexes
$R_{T\cup\{v_0\}}$. By \eqref{eq:RT} and $|T|\le a$,
\[
L_\kk(B)
\le\max_{1\le r\le a}\bigl(r(d-r)+r-1\bigr)
=\max_{1\le r\le a}r(d+1-r)-1
\le\left\lfloor\frac{(d+1)^2}{4}\right\rfloor-1.
\]
Since $\lfloor(d+3)^2/4\rfloor=\lfloor(d+1)^2/4\rfloor+d+2$ and $a\le d$,
every $j\ge D$ satisfies
\[
j-a-2\;\ge\;\left\lfloor\frac{(d+1)^2}{4}\right\rfloor-1\;\ge\;L_\kk(B),
\]
and the middle term is nonnegative because $d\ge1$. Thus
$\Ht_{j-a-2}(B;\kk)=0$, and all the relative homology groups above vanish
for $j\ge D$. The induction is complete.
\end{proof}

\section{Sums and joins}\label{sec:joins}

Sums of ideals in disjoint sets of variables correspond to
joins of complexes. In this section, we recall how the squarefree
symbolic powers of such sums are expressed in terms of those of their
summands and how their regularity is bounded, and show that
Conjecture~\ref{conj:leray} is preserved under joins.

Throughout, $R_1$ and $R_2$ are polynomial rings over $\kk$
in disjoint sets of variables, $S=R_1\otimes_\kk R_2$, and
$I\subsetneq R_1$, $J\subsetneq R_2$ are squarefree monomial
ideals. Write $I=I_{K_1}$ and $J=I_{K_2}$, where
$K_1$ and $K_2$ are nonvoid complexes on disjoint sets.
Extending $I$ and $J$ to $S$, and denoting the extensions
by the same letters, we have $I+J=I_{K_1*K_2}$.
Put $\sym{I}{0}=R_1$ and $\sym{J}{0}=R_2$, and for
$a,b\ge1$ set
$f_a=\reg(R_1/\sym{I}{a})$ and
$g_b=\reg(R_2/\sym{J}{b})$. Note that $f_a=0$ if $a>\hgt I$ and $g_b=0$ if $b>\hgt J$, since the corresponding squarefree symbolic powers vanish.%

\begin{theorem}[\cite{HNTT20,CDRS}]\label{thm:sum}
For every integer $k\ge1$, the following hold.
\begin{enumerate}[label=\textup{(\arabic*)}]
\item $\sym{(I+J)}{k}
=\sum_{a+b=k}\sym{I}{a}\,\sym{J}{b}$, where $a,b$ range over integers with $a,b\ge0$.
\item The following bound holds:
\[
\reg(S/\sym{(I+J)}{k})
\le\max\bigl(
\{f_a+g_b+1:a+b=k\}
\cup\{f_a+g_b:a+b=k+1\}\bigr),
\]
where $a,b$ range over positive integers.
\end{enumerate}
\end{theorem}

Part~(1) is \cite[Definition~3.1(d) and Proposition~4.2]{CDRS}; it is the
squarefree part of the binomial expansion
$(I+J)^{(k)}=\sum_{a+b=k}I^{(a)}J^{(b)}$ of
\cite[Theorem~3.4]{HNTT20}. Part~(2) follows from
\cite[Theorems~2.6 and~2.7]{CDRS}, since squarefree symbolic powers form
Tor-vanishing filtrations \cite[Corollary~2.10]{CDRS}.
As for symbolic powers \cite[Theorem~5.11]{HNTT20}, the bound in Theorem~\ref{thm:sum}(2)
is an
equality \cite[Theorem~2.7]{CDRS} for $1\le k\le\hgt I+\hgt J$ if $a$ and $b$ are
restricted to
$1\le a\le\hgt I$ and $1\le b\le\hgt J$; we use only the inequality.

To apply Theorem~\ref{thm:sum} to the
conjectured bound, we need the following
inequalities for Erd\H{o}s--Gallai numbers.

\begin{lemma}\label{lem:etamerge}
For all integers $r_1,r_2,a,b\ge1$,
$\eta(\max\{r_1,r_2\},a+b)
\ge\eta(r_1,a)+\eta(r_2,b)
$
and
$\eta(r_1+r_2-1,a+b-1)
\ge\eta(r_1,a)+\eta(r_2,b)-1.
$
\end{lemma}

\begin{proof}
Let $\cH_1$ be an $r_1$-uniform
$a$-critical hypergraph with
$\eta(r_1,a)$ vertices, and let
$\cH_2$ be an $r_2$-uniform
$b$-critical hypergraph with
$\eta(r_2,b)$ vertices, on disjoint
vertex sets. 
Their disjoint union is $(a+b)$-critical,
and its edges have cardinality at most
$\max\{r_1,r_2\}$. The first inequality
follows from Lemma~\ref{lem:eta}(2).

For the second, fix an edge $e$ of
$\cH_1$ and a vertex $y$ of $\cH_2$.
Let $E_y$ be the nonempty family of
edges of $\cH_2$ containing $y$.
Construct $\cH$ by substituting the
edge $e$ for $y$ in $\cH_2$ and
adjoining the remaining edges of
$\cH_1$. Thus its edges are those
of $\cH_1-e$, those of $\cH_2$
outside $E_y$, and
$e\cup(h\setminus\{y\})$ for
$h\in E_y$. Every edge has cardinality
at most $r_1+r_2-1$. Moreover, $\cH$
has no isolated vertices and has
$\eta(r_1,a)+\eta(r_2,b)-1$ vertices.
By criticality, $\cH_1-e$ has a
transversal of size $a-1$, and every
such transversal misses $e$.
Similarly, for every edge $h$ of
$\cH_2$, the hypergraph $\cH_2-h$
has a transversal of size $b-1$,
and every such transversal misses $h$.
In particular, there is a transversal
$Y$ of $\cH_2$ of size $b$ containing
$y$: choose $h\in E_y$ and add $y$
to a transversal of $\cH_2-h$ of
size $b-1$. Then $Y\setminus\{y\}$
meets every edge of $\cH_2$
outside $E_y$.
Let $X=X_1\sqcup X_2$ be a
transversal of $\cH$, with $X_i$
contained in the vertex set of
$\cH_i$. If $X_1$ meets $e$, then
$X_1$ is a transversal of $\cH_1$,
while $X_2\cup\{y\}$ is a transversal
of $\cH_2$. Hence $|X_1|\ge a$
and $|X_2|\ge b-1$.
If $X_1$ misses $e$, then $X_1$
is a transversal of $\cH_1-e$
and $X_2$ is a transversal of
$\cH_2$, so $|X_1|\ge a-1$
and $|X_2|\ge b$.
Thus $\tau(\cH)\ge a+b-1$.
Conversely, the union of a transversal
of $\cH_1$ of size $a$ with
$Y\setminus\{y\}$ is a transversal
of $\cH$. Therefore
$\tau(\cH)=a+b-1$.

It remains to show that deleting any
edge of $\cH$ decreases its transversal
number. We distinguish three cases.

\begin{enumerate}[label=\textup{(\roman*)}]
\item Let $e'\ne e$ be an edge of
$\cH_1$. A transversal of
$\cH_1-e'$ of size $a-1$ meets $e$.
Its union with $Y\setminus\{y\}$
is a transversal of $\cH-e'$
of size $a+b-2$.

\item Let $e'$ be an edge of
$\cH_2$ outside $E_y$, and let
$X_2$ be a transversal of
$\cH_2-e'$ of size $b-1$.
If $y\notin X_2$, take the union
of $X_2$ with a transversal of
$\cH_1-e$ of size $a-1$.
If $y\in X_2$, take the union
of $X_2\setminus\{y\}$ with a
transversal of $\cH_1$ of size $a$.
In either case, the resulting set
is a transversal of $\cH-e'$
of size $a+b-2$.

\item Let $e'=e\cup(h\setminus\{y\})$
for some $h\in E_y$.
Choose a transversal $X_2$ of
$\cH_2-h$ of size $b-1$.
Since $X_2$ misses $h$, it does
not contain $y$. Its union with
a transversal of $\cH_1-e$ of
size $a-1$ meets every edge of
$\cH-e'$ and has size $a+b-2$.
\end{enumerate}

Thus $\cH$ is $(a+b-1)$-critical,
and the second inequality follows
from Lemma~\ref{lem:eta}(2).
\end{proof}

\begin{corollary}\label{cor:joins}
Let $d_1,d_2\ge0$ and $d=d_1+d_2$.
Suppose that $\reg(R_1/I)\le d_1$ and
$\reg(R_2/J)\le d_2$, and that the
conclusion of Conjecture~\ref{conj:reg}
holds for $I$ with $d_1$ and for $J$
with $d_2$. Then $\reg(S/(I+J))\le d$,
and the conclusion of
Conjecture~\ref{conj:reg} holds for
$I+J$ with $d$.
In particular, if $K_1$ and $K_2$
satisfy Conjecture~\ref{conj:leray}
over $\kk$, then so does $K_1*K_2$.
\end{corollary}

\begin{proof}
Since
$S/(I+J)\cong(R_1/I)\otimes_\kk(R_2/J)$,
regularity is additive:
\begin{equation}\label{eq:regjoin}
\reg(S/(I+J))
=\reg(R_1/I)+\reg(R_2/J)\le d.
\end{equation}
Put $r_i=d_i+1$, so that
$r_1+r_2-1=d+1$.
By hypothesis,
$f_a\le\eta(r_1,a)-1$ and
$g_b\le\eta(r_2,b)-1$
for all $a,b\ge1$.
Let $k\ge1$.
If $a+b=k$, then
Lemmas~\ref{lem:etamerge}
and~\ref{lem:eta}(1) give
\[
\begin{aligned}
f_a+g_b+1
\le\eta(r_1,a)+\eta(r_2,b)-1
\le\eta(\max\{r_1,r_2\},k)-1
\le\eta(d+1,k)-1.
\end{aligned}
\]
If $a+b=k+1$, then
Lemma~\ref{lem:etamerge} gives
\[
\begin{aligned}
f_a+g_b
\le\eta(r_1,a)+\eta(r_2,b)-2
\le\eta(d+1,k)-1.
\end{aligned}
\]
Theorem~\ref{thm:sum}(2) therefore
yields
$\reg(S/\sym{(I+J)}{k})
\le\eta(d+1,k)-1$.

For the last assertion, take
$d_i=L_\kk(K_i)$.
By \eqref{eq:regjoin},
$L_\kk(K_1*K_2)=d$.
Since $\eta(r,k)$ is nondecreasing
in $r$ by Lemma~\ref{lem:eta}(1),
it suffices to verify
Conjecture~\ref{conj:leray} with
$d$ equal to the Leray number.
The first assertion applies.
\end{proof}

\begin{remark}\label{rem:joincoll}
We do not know whether Conjecture~\ref{conj:coll} is preserved under joins,
although its hypothesis is preserved, since
$C(K_1*K_2)\le C(K_1)+C(K_2).$ Indeed, let $K_1=X_0,X_1,\dots,X_p=\emptyset$ and $K_2=Y_0,Y_1,\dots,Y_q=\emptyset$ be collapsing sequences, where $X_i$ arises from $X_{i-1}$ by an elementary collapse at $\sigma_i$ with unique maximal face $M_i$, and $Y_j$ arises from $Y_{j-1}$ by a collapse at $\alpha_j$ with unique maximal face $N_j$. In $K_1*K_2$, remove the faces containing $\sigma_i\cup\alpha_j$ for the pairs $(i,j)$ in lexicographic order. Every face of $K_1$ contains some $\sigma_{i'}$, and every face of $K_2$ some $\alpha_{j'}$. Hence a face $\rho\cup\beta\supseteq\sigma_i\cup\alpha_j$ that survives the earlier steps satisfies $\rho\in X_{i-1}$ and $\beta\in Y_{j-1}$, so $\rho\subseteq M_i$ and $\beta\subseteq N_j$; conversely, every face of $M_i\cup N_j$ containing $\sigma_i\cup\alpha_j$ survives. Thus $\sigma_i\cup\alpha_j$ lies in the unique maximal face $M_i\cup N_j$ of the current complex, each step is an elementary collapse at a face of cardinality at most $C(K_1)+C(K_2)$, and every face of $K_1*K_2$ is eventually removed.

The $\theta$-number of~\cite{BC22} is additive under joins
\cite[Theorem~11]{BC22}, but whether equality always holds in the
inequality above remains open, as also noted in \cite{BC22}. The obstacle to a collapsibility analogue of
Corollary~\ref{cor:joins} is that
$\T_t(K_1*K_2)=\bigcup_{a+b=t}\T_a(K_1)*\T_b(K_2)$
is a union of joins, for which we have no analogue of the filtration
argument behind Theorem~\ref{thm:sum}.
\end{remark}

\section{Top-dimensional homology and matroids}
\label{sec:top}

By Remark~\ref{rem:dimbound}, every nonvoid complex
$K$ satisfies
$L_\kk(\T_t(K))\le C(\T_t(K))
\le\dim K+1+t$ $(t\ge0).$
In this section, we show that these bounds are
attained whenever $\Ht_{\dim K}(K;\kk)\ne0$
and $t<\hgt I_K$. The proof constructs
top-dimensional cycles explicitly.
As an application, we determine the invariants
of tolerance complexes of matroids and the
regularity of squarefree symbolic powers of
circuit ideals.

Throughout, $\widetilde C_\bullet$ denotes the
augmented oriented simplicial chain complex
over $\kk$. A nonzero cycle in the top dimension
of a complex represents a nonzero homology
class, since there are no boundaries in that
dimension.

\begin{lemma}\label{lem:topLeray}
Let $X$ be a nonvoid complex. Then
$L_\kk(X)\le\dim X+1$, with equality if
and only if $\Ht_{\dim X}(X;\kk)\ne0$.
Equivalently,
$\reg(S/I_X)\le\dim(S/I_X)$, with equality
if and only if $\Ht_{\dim X}(X;\kk)\ne0$.
\end{lemma}

\begin{proof}
Put $D=\dim X$; the case $D=-1$, where $X=\{\emptyset\}$, is immediate.
Every induced subcomplex $X[U]$ has dimension at most $D$, so
$L_\kk(X)\le D+1$. Moreover, $\Ht_D(X[U];\kk)$ is the space of $D$-cycles of
$X[U]$, and each such cycle is also a $D$-cycle of $X$; hence
$L_\kk(X)=D+1$ if and only if $\Ht_D(X;\kk)\ne0$. The algebraic statement
follows from Hochster's formula and $\dim(S/I_X)=D+1$.
\end{proof}

To construct cycles in tolerance complexes,
we use the contracting homotopy of a cone.
For $v\in W$, define
$c_v\colon\widetilde C_i(2^W)
\to\widetilde C_{i+1}(2^W)$ by
\[
c_v[f]=
\begin{cases}
[v,f],&v\notin f,\\
0,&v\in f.
\end{cases}
\]
Here $[v,f]$ denotes the oriented simplex obtained by placing $v$ first, and $c_v[\emptyset]=[v]$ in degree $-1$. Then
\begin{equation}\label{eq:cone}
\partial c_v+c_v\partial=\mathrm{id}
\qquad\text{on }\widetilde C_\bullet(2^W).
\end{equation}
Indeed, if $v\notin f$, then
$\partial[v,f]=[f]-[v,\partial f]$
and $c_v\partial[f]=[v,\partial f]$.
If $f=[v,g]$, then $c_v[f]=0$ and
$c_v\partial[f]=[f]$.

\begin{lemma}\label{lem:onestep}
Let $X$ be a nonvoid complex on $W$ with
$\Ht_{\dim X}(X;\kk)\ne0$ and
$\hgt I_X\ge2$. Then
$\Ht_{\dim X+1}(\T_1(X);\kk)\ne0$.
\end{lemma}

\begin{proof}
Put $D=\dim X$. By
Proposition~\ref{prop:height},
$|W|=D+1+\hgt I_X\ge D+3$ and
$\dim\T_1(X)=D+1$.
It therefore suffices to construct
a nonzero $(D+1)$-cycle in $\T_1(X)$.
Let $z\in\widetilde C_D(X)$ be a
nonzero cycle. For every $v\in W$,
the chain $c_vz$ is supported on
faces of $\T_1(X)$, and
$\partial(c_vz)=z$ by \eqref{eq:cone}.
Thus $c_uz-c_vz$ is a $(D+1)$-cycle
in $\T_1(X)$ for all $u,v\in W$.

Suppose all these cycles vanish, and
let $c$ be the common value of the
chains $c_vz$. Every simplex in the
support of $c_vz$ contains $v$.
Hence every simplex in the support
of $c$ contains all vertices of $W$.
But these simplices have cardinality
$D+2<|W|$, so $c=0$.
This contradicts $\partial c=z\ne0$.
\end{proof}

\begin{theorem}\label{thm:top}
Let $K$ be a nonvoid complex on $W$,
put $D=\dim K$, and suppose that
$\Ht_D(K;\kk)\ne0$. Then, for every
integer $0\le t<\hgt I_K$,
$\Ht_{D+t}(\T_t(K);\kk)\ne0$.
Consequently,
\[
C(\T_t(K))=L_\kk(\T_t(K))
=\dim\T_t(K)+1=L_\kk(K)+t.
\]
\end{theorem}

\begin{proof}
We induct on $t$. The case $t=0$
is the hypothesis. Let
$1\le t<\hgt I_K$ and put
$X=\T_{t-1}(K)$.
By Proposition~\ref{prop:height},
$\dim X=D+t-1$ and
$\hgt I_X=\hgt I_K-t+1\ge2$.
By induction,
$\Ht_{D+t-1}(X;\kk)\ne0$.
Lemma~\ref{lem:onestep} therefore gives
$\Ht_{D+t}(\T_1(X);\kk)\ne0$.
Since $\T_1(X)=\T_t(K)$, the first
assertion follows.
By Lemma~\ref{lem:topLeray} and
Proposition~\ref{prop:height},
$L_\kk(\T_t(K))=\dim\T_t(K)+1=D+t+1$.
Moreover,
$L_\kk(\T_t(K))
\le C(\T_t(K))
\le\dim\T_t(K)+1$
by Lemma~\ref{lem:hLC} and
\eqref{eq:dimC}, so equality holds
throughout. Finally,
$L_\kk(K)=D+1$ by
Lemma~\ref{lem:topLeray}.
\end{proof}

The restriction $t<\hgt I_K$ is
necessary, since $\T_t(K)=2^W$
for $t\ge\hgt I_K$.
The hypothesis of Theorem~\ref{thm:top}
holds, for example, for Gorenstein$^*$
complexes over $\kk$, closed orientable
pseudomanifolds, Cohen--Macaulay
complexes whose $h$-vector has nonzero
last entry, and independence complexes
of matroids without coloops.

We conclude with matroid complexes.
Let $M$ be a matroid on a finite set
$E$ with rank function $r_M$. Write
$r(M)=r_M(E)$,
$\nu_M(U)=|U|-r_M(U)$ for $U\subseteq E$,
and $\nu(M)=\nu_M(E)$.
Let $c(M)$ denote the number of
coloops of $M$. The independence
complex $\Ind(M)$ has
\emph{circuit ideal}
$I_M=(x^C:C\text{ is a circuit of }M)$
as its Stanley--Reisner ideal, with
$\hgt I_M=\nu(M)$.
For $t\ge0$, the \emph{elongation}
$M^{+t}$ is the matroid on $E$
whose independent sets are the sets
$X\subseteq E$ satisfying
$\nu_M(X)\le t$ \cite{Oxl11}.
Its rank function is
$r_{M^{+t}}(U)
=\min\{|U|,r_M(U)+t\}$.

As observed by Kim and Lew~\cite{KL24}, and in terms of squarefree symbolic
powers by Mantero and Nguyen~\cite[Corollary~3.21]{MN24}, tolerance complexes
of matroid complexes are matroid complexes:
$\T_t(\Ind(M))=\Ind(M^{+t})$, since a set $X$ belongs to $\T_t(\Ind(M))$
exactly when $r_M(X)\ge|X|-t$, that is, $\nu_M(X)\le t$.

By Proposition~\ref{prop:sym},
$\sym{I_M}{t+1}=I_{M^{+t}}$.
Thus every squarefree symbolic power of a circuit
ideal is again a circuit ideal, and its quotient
ring is Cohen--Macaulay \cite{Bjo92}.

\begin{theorem}[{cf.\ \cite{KL24,MN24}, \cite[Theorem~30]{BC22}}]\label{thm:matroid}
Let $K=\Ind(M)$ and $0\le t<\nu(M)$. Then
\[
C(\T_t(K))=L_\kk(\T_t(K))=r(M)-c(M)+t=C(K)+t.
\]
In algebraic terms, for $1\le k\le\nu(M)$,
$\reg(S/\sym{I_M}{k})=r(M)-c(M)+k-1 .$
Moreover, $S/\sym{I_M}{k}$ is Cohen--Macaulay with
$\operatorname{pd}(S/\sym{I_M}{k})=\nu(M)-k+1$ and $a$-invariant $-c(M)$.
In particular, Conjectures~\ref{conj:leray} and~\ref{conj:coll} hold for
matroid complexes.
\end{theorem}

\begin{proof}
Let $E_0$ be the set of coloops of $M$ and put $M'=M|(E\setminus E_0)$.
Then $K=2^{E_0}*\Ind(M')$, $r(M')=r(M)-c(M)$, and $\nu(M')=\nu(M)$.
For every complex $X$ on $E\setminus E_0$, we have
$\T_t(2^{E_0}*X)=2^{E_0}*\T_t(X)$; indeed, a witnessing face of
$2^{E_0}*X$ inside $\sigma$ may be enlarged to contain $\sigma\cap E_0$.
Moreover, $C(2^{E_0}*X)=C(X)$ by \cite[Lemma~3.3]{KL22}, while
$L_\kk(2^{E_0}*X)=L_\kk(X)$ by~\eqref{eq:regjoin}.
Thus, for the first formula, we may replace $M$ by $M'$ and assume that
$M$ has no coloops. 
In this case, $\dim K=r(M)-1$ and
$\dim_\kk\Ht_{r(M)-1}(K;\kk)=T_M(0,1)>0$, where $T_M$ is the Tutte
polynomial \cite[Section~7.8]{Bjo92}; this remains valid when $r(M)=0$,
using reduced homology in degree $-1$.
Since $t<\nu(M)=\hgt I_M$, Theorem~\ref{thm:top} gives
$C(\T_t(K))=L_\kk(\T_t(K))=r(M)+t$.
Taking $t=0$ also gives $C(K)=r(M)$.
Restoring the coloops proves the first formula.

For the algebraic statements, put $t=k-1$. The regularity formula follows
from the first formula and Proposition~\ref{prop:sym}, since
$L_\kk(\T_t(K))=\reg(S/\sym{I_M}{k})$.
By Proposition~\ref{prop:sym} and the identity $\T_t(\Ind(M))=\Ind(M^{+t})$,
the ideal $\sym{I_M}{k}$ is the circuit ideal of the matroid $M^{+t}$, so
$S/\sym{I_M}{k}$ is Cohen--Macaulay \cite{Bjo92}. Hence its projective
dimension equals the height of $\sym{I_M}{k}$, which is $\nu(M)-k+1$ by
Proposition~\ref{prop:height}. Its Krull dimension is $r(M)+k-1$, so its
$a$-invariant is $(r(M)-c(M)+k-1)-(r(M)+k-1)=-c(M)$.
Finally, suppose that $K$ is $d$-Leray; in particular, this holds if $K$ is
$d$-collapsible. Then $L_\kk(K)=r(M)-c(M)\le d$; if $\nu(M)=0$, then
$K=2^E$ and both sides of this equality are $0$.
For $t<\nu(M)$, Lemma~\ref{lem:eta}(1) gives
\[
C(\T_t(K))=r(M)-c(M)+t\le d+t\le(d+1)(t+1)-1\le\eta(d+1,t+1)-1 .
\]
For $t\ge\nu(M)$, we have $\T_t(K)=2^E$. This completes the proof.
\end{proof}

Theorem~\ref{thm:matroid} can also be derived from known results on matroid complexes
\cite{KL24,MN24,Bjo92} and \cite[Theorem~30]{BC22}; see also \cite{DFKT24}; the proof above obtains it directly
from
Theorem~\ref{thm:top}.

\section{Stable and shifted complexes}\label{sec:shifted}

In this section, we study complexes whose Stanley--Reisner
ideals are squarefree stable. For these complexes, tolerance
satisfies linear bounds on the three invariants, while
for shifted complexes we obtain an exact description.
Throughout, $W=[n]$ and $S=\kk[x_1,\dots,x_n]$.
For a monomial $u\ne1$, let $\min(u)$ and $\max(u)$
denote the least and largest indices of variables
dividing $u$, respectively. We write $G(I)$ for the
minimal monomial generating set of a monomial ideal $I$. 
Following Aramova, Herzog, and Hibi \cite{AHH98},
a squarefree monomial ideal $I$ is
\emph{squarefree stable} if
$x_j(u/x_{\max(u)})\in I$ for every squarefree
monomial $u\in I$ and every $j<\max(u)$ with
$x_j\nmid u$. It suffices to check this condition
for $u\in G(I)$. The ideal $I$ is
\emph{squarefree strongly stable} if
$x_j(u/x_i)\in I$ whenever $u\in I$ is squarefree,
$x_i\mid u$, $x_j\nmid u$, and $j<i$.
Every squarefree strongly stable ideal is
squarefree stable, but not conversely.

A nonvoid complex $K$ on $[n]$ is called
\emph{stable} if $I_K$ is squarefree stable.
Equivalently, if $\sigma$ is a nonface of $K$,
$j<\max\sigma$, and $j\notin\sigma$, then
$(\sigma\setminus\{\max\sigma\})\cup\{j\}$
is also a nonface. A complex $K$ is
\emph{shifted} if $f\in K$, $x\in f$, $y<x$,
and $y\notin f$ imply
$(f\setminus\{x\})\cup\{y\}\in K$.
After reversing the order of $[n]$, this
is equivalent to $I_K$ being squarefree
strongly stable. Thus shifted complexes
become stable after reversing the vertex
order, and the results below for stable
complexes apply to them.

\begin{lemma}\label{lem:stablemove}
Let $K$ be stable, $f\in K$, $y\in f$,
and $M\in[n]$ with $M>\max f$. Then
$(f\setminus\{y\})\cup\{M\}\in K$.
\end{lemma}

\begin{proof}
Put $\rho=(f\setminus\{y\})\cup\{M\}$.
Then $\max\rho=M$, $y<M$, and
$y\notin\rho$. If $\rho$ were a nonface,
stability would imply that
$(\rho\setminus\{M\})\cup\{y\}=f$
is a nonface, a contradiction.
\end{proof}

\begin{lemma}\label{lem:Tstable}
If $K$ is stable, then $\T_t(K)$ is stable
for every $t\ge0$.
\end{lemma}

\begin{proof}
Since $K$ is nonvoid, every set of cardinality at most $t$
belongs to $\T_t(K)$. A set $\sigma$
with $|\sigma|>t$ is a nonface of
$\T_t(K)$ if and only if every subset
of $\sigma$ of cardinality $|\sigma|-t$
is a nonface of $K$.
Let $\sigma$ be a nonface of $\T_t(K)$,
put $m=\max\sigma$, and let $j<m$
with $j\notin\sigma$. Set
$\sigma'=(\sigma\setminus\{m\})\cup\{j\}$.
Let $\rho'\subseteq\sigma'$ have
cardinality $|\sigma|-t$.
If $j\notin\rho'$, then
$\rho'\subseteq\sigma$, so $\rho'$
is a nonface of $K$.
Otherwise, write $\rho'=\rho\cup\{j\}$
with $\rho\subseteq\sigma\setminus\{m\}$.
Then $\rho\cup\{m\}$ is a nonface of
$K$ with maximum $m$, so stability
implies that $\rho'$ is a nonface.
Thus $\sigma'\notin\T_t(K)$.
\end{proof}

\begin{proposition}\label{prop:stableCLh}
If $X$ is a nonvoid stable complex, then
$C(X)=L_\kk(X)=h(X)$ for every field $\kk$.
\end{proposition}

The equality $L_\kk(X)=h(X)$ is known. If $I_X=0$, then $X=2^{[n]}$ and both
sides vanish. Otherwise, it follows from the fact that the regularity of a
squarefree stable ideal equals the largest degree of a minimal generator
\cite{AHH98}; see also \cite[Example~1.2 and Corollary~1.6]{HT02}.
The new content of Proposition~\ref{prop:stableCLh} is the equality $C(X)=h(X)$.%

\begin{proof}
By Lemma~\ref{lem:hLC}, it suffices to show that $C(X)\le h(X)$.
We induct on $n$. Since stability depends only on the order of the ground
set, a stable complex on $\{2,\dots,n\}$ becomes a stable complex on $[n-1]$
under the order-preserving relabelling, so the induction hypothesis applies
to it. 
If $h(X)=0$, then $X$ has no missing face of cardinality at least two, so $X$
is a simplex on its vertex set and $C(X)=0$. Suppose $d=h(X)\ge1$.
If $1$ is not a vertex of $X$, then $X$, regarded as a complex on
$\{2,\dots,n\}$, is stable with the same invariants, and the assertion
follows by induction. Now suppose that $1$ is a vertex of $X$.
The deletion $X\setminus1=X[\{2,\dots,n\}]$ is a nonvoid stable complex on
$\{2,\dots,n\}$. Indeed, its nonfaces are the nonfaces of $X$ contained in
$\{2,\dots,n\}$, and stability of $X$ applied with $j\ge2$ preserves this
property. Its missing faces are missing faces of $X$, so $h(X\setminus1)\le d$,
and induction gives $C(X\setminus1)\le d$. 
The link $\lk(X,1)$ is a nonvoid complex on $\{2,\dots,n\}$, and all its
nonfaces are nonempty, since $1$ is a vertex of $X$. It is stable. Indeed, let
$\rho$ be a nonface of $\lk(X,1)$, that is, $\rho\cup\{1\}\notin X$, and let
$j<\max\rho$ with $j\notin\rho\cup\{1\}$. Since $\max(\rho\cup\{1\})=\max\rho$,
stability of $X$ shows that $(\rho\setminus\{\max\rho\})\cup\{j,1\}\notin X$;
that is, $(\rho\setminus\{\max\rho\})\cup\{j\}$ is a nonface of $\lk(X,1)$.

Let $\rho$ be a missing face of $\lk(X,1)$, and put $M=\max\rho$; thus
$\rho\ne\emptyset$ and $M\ge2$. Since $\rho\setminus\{M\}\in\lk(X,1)$, we have
$(\rho\setminus\{M\})\cup\{1\}\in X$. If $\rho\notin X$, then stability
applied to the nonface $\rho$ with $j=1$ would give
$(\rho\setminus\{M\})\cup\{1\}\notin X$, a contradiction. Thus $\rho\in X$.
Every proper subset of $\rho\cup\{1\}$ is either contained in $\rho$ or of the
form $\rho'\cup\{1\}$ with $\rho'\subsetneq\rho$, and hence belongs to $X$,
while $\rho\cup\{1\}\notin X$. Hence $\rho\cup\{1\}$ is a missing face of $X$,
so $|\rho|\le d$. Consequently, $h(\lk(X,1))\le d-1$, and induction gives
$C(\lk(X,1))\le d-1$. Since $1\in V(X)$, Lemma~\ref{lem:tancer} yields
$C(X)\le\max\{C(X\setminus1),\,C(\lk(X,1))+1\}\le d.$
\end{proof}

\begin{theorem}\label{thm:stable}
Let $K$ be a nonvoid stable complex
and let $t\ge0$. Then
\[
C(\T_t(K))=L_\kk(\T_t(K))
=h(\T_t(K))\le h(K)+t=C(K)+t.
\]
In particular, Conjectures~\ref{conj:leray}
and~\ref{conj:coll} hold for
Stanley--Reisner complexes of
squarefree stable ideals.
\end{theorem}

\begin{proof}
The equalities follow from
Lemma~\ref{lem:Tstable} and
Proposition~\ref{prop:stableCLh}.
If $\T_t(K)=2^{[n]}$, the inequality
is immediate. Otherwise, let $\sigma$
be a missing face of $\T_t(K)$, and put
$s=|\sigma|-t\ge1$ and $M=\max\sigma$.
Since $\sigma\setminus\{M\}\in\T_t(K)$,
it contains a face $f\in K$ of
cardinality $s-1$.
The set $f\cup\{M\}\subseteq\sigma$
has cardinality $s$, so it is a
nonface of $K$. Its proper subsets
avoiding $M$ are contained in $f$,
while those containing $M$ are
contained in
$(f\setminus\{y\})\cup\{M\}$
for some $y\in f$, which belongs
to $K$ by Lemma~\ref{lem:stablemove}.
Thus $f\cup\{M\}$ is a missing face
of $K$, and $s\le h(K)+1$.
It follows that
$|\sigma|-1\le h(K)+t$.
Finally, if $K$ is $d$-Leray or
$d$-collapsible, then
$C(K)=h(K)\le d$. Hence
$C(\T_t(K))\le d+t
\le(d+1)(t+1)-1
\le\eta(d+1,t+1)-1$
by Lemma~\ref{lem:eta}(1).
\end{proof}

The inequality is sharp. Indeed, if $K$
consists of all subsets of $[n]$ of
cardinality at most $d$, with
$n\ge d+t+1$, then
$h(K)=d$ and $h(\T_t(K))=d+t$.
Algebraically, if $I\ne0$ is squarefree
stable and $\sym{I}{k}\ne0$, then
$\sym{I}{k}$ is squarefree stable and
$\reg\sym{I}{k}\le\reg I+k-1$,
by Lemma \ref{lem:Tstable}, Theorem \ref{thm:stable} and Proposition~\ref{prop:sym}.
For shifted complexes, we can describe
the tolerance complexes, and hence
squarefree symbolic powers of squarefree
strongly stable ideals, explicitly.
For $f=\{f_1<\dots<f_q\}$ and
$g=\{g_1<\dots<g_q\}$, write
$g\le f$ if $g_i\le f_i$ for all $i$.
It is well known that if $K$ is shifted,
$f\in K$, and $g\le f$, then $g\in K$;
see \cite{Kal02}.

\begin{lemma}\label{lem:Tshift}
If $K$ is shifted, then so is
$\T_t(K)$ for every $t\ge0$.
\end{lemma}

\begin{proof}
Let $f\in K$, $f\subseteq\sigma$,
and $|\sigma\setminus f|\le t$.
Let $x\in\sigma$, $y<x$, and
$y\notin\sigma$, and put
$\sigma'=(\sigma\setminus\{x\})\cup\{y\}$.
If $x\notin f$, put $f'=f$;
otherwise, put
$f'=(f\setminus\{x\})\cup\{y\}$.
In both cases, $f'\in K$,
$f'\subseteq\sigma'$, and
$|\sigma'\setminus f'|
=|\sigma\setminus f|\le t$.
Thus $\sigma'\in\T_t(K)$.
\end{proof}

\begin{theorem}\label{thm:shiftexact}
Let $K$ be a nonvoid shifted complex
on $[n]$, and let $t\ge0$.
The missing faces of $\T_t(K)$ are
exactly the sets $B\sqcup R$, where
$B$ is a missing face of $K$ and
$R$ is a $t$-element subset of
$\{\max B+1,\dots,n\}$.
Distinct pairs $(B,R)$ give distinct sets.
If no missing face $B$ of $K$ satisfies
$\max B\le n-t$, then
$\T_t(K)=2^{[n]}$.
Otherwise, for every field $\kk$,
\[
\begin{aligned}
C(\T_t(K))&=L_\kk(\T_t(K))
=h(\T_t(K))\\
&=t-1+
\max\{|B|:B\text{ is a missing face of }K,
\ \max B\le n-t\}.
\end{aligned}
\]
\end{theorem}

\begin{proof}
Let $\sigma=\{a_1<\dots<a_m\}$ be
a missing face of $\T_t(K)$.
Since every set of cardinality at
most $t$ belongs to $\T_t(K)$,
we have $m\ge t+1$.
Put $B=\{a_1,\dots,a_{m-t}\}$.
Then $B\notin K$, since otherwise
$\sigma\in\T_t(K)$.
As $\sigma\setminus\{a_1\}\in\T_t(K)$,
it contains a face $f\in K$ with
$|f|=m-t-1$.
Since $f\subseteq\{a_2,\dots,a_m\}$,
we have $B\setminus\{a_1\}\le f$,
and hence $B\setminus\{a_1\}\in K$.
For every $j$,
$B\setminus\{a_j\}
\le B\setminus\{a_1\}$,
so $B\setminus\{a_j\}\in K$.
Thus $B$ is a missing face of $K$,
and $R=\sigma\setminus B$ consists
of $t$ elements larger than $\max B$.

Conversely, let $B$ be a missing
face of $K$, and let
$R\subseteq\{\max B+1,\dots,n\}$
with $|R|=t$. Put $\sigma=B\sqcup R$.
The increasing enumeration of $\sigma$
begins with that of $B$, so every
$|B|$-subset $f$ of $\sigma$ satisfies
$B\le f$, and hence $f\notin K$.
Thus $\sigma\notin\T_t(K)$.
For any $x\in\sigma$, the set
$\sigma\setminus\{x\}$ contains
$B\setminus\{b\}$, where $b=x$
if $x\in B$, and $b$ is any element
of $B$ otherwise.
This is a face of $K$ and misses
exactly $t$ elements of
$\sigma\setminus\{x\}$.
Hence $\sigma\setminus\{x\}\in\T_t(K)$,
so $\sigma$ is a missing face.
Given $\sigma$ and $t$, the set $B$
consists of the $|\sigma|-t$ smallest
elements of $\sigma$, so the pair
$(B,R)$ is uniquely determined.
The formula for the invariants
follows from Lemma~\ref{lem:Tshift}
and Proposition~\ref{prop:stableCLh},
after reversing the order of $[n]$.
\end{proof}

In algebraic terms,
Theorem~\ref{thm:shiftexact}
determines the minimal generators,
graded Betti numbers, regularity,
and projective dimension of all
squarefree symbolic powers of
squarefree strongly stable ideals.

\begin{corollary}\label{cor:stable}
Let $I\subsetneq S$ be a nonzero
squarefree strongly stable ideal,
and let $k\ge1$. Then $\sym{I}{k}$
is squarefree strongly stable, and
\[
G(\sym{I}{k})
=\{uw:u\in G(I),\ \deg w=k-1,\
\max(w)<\min(u)\},
\]
where $w$ ranges over squarefree
monomials, with $\max(1)=0$,
and distinct pairs $(u,w)$ give
distinct generators.
Consequently, for all $i,j$,
$$\beta_{i,i+j}(\sym{I}{k})
=\sum_{\substack{u\in G(I)\\
\deg u=j-k+1}}
\binom{\min(u)-1}{k-1}
\binom{\max(u)-j}{i}.
$$
If $\sym{I}{k}\ne0$, then
\[
\begin{aligned}
\reg\sym{I}{k}
&=\max\{\deg u:u\in G(I),\
\min(u)\ge k\}+k-1,\\
\operatorname{pd}\sym{I}{k}
&=\max\{\max(u)-\deg u:
u\in G(I),\ \min(u)\ge k\}-k+1.
\end{aligned}
\]
In particular,
$\reg\sym{I}{k}\le\reg I+k-1$,
with equality for squarefree Veronese
ideals whenever $\sym{I}{k}\ne0$.
\end{corollary}

\begin{proof}
Reversing the order of the variables
turns $I$ into the Stanley--Reisner
ideal of a shifted complex $K$,
with $\max B$ corresponding to
$\min(u)$. Thus
Theorem~\ref{thm:shiftexact}
and Proposition~\ref{prop:sym}
give the minimal generators.
Moreover, $\sym{I}{k}$ is
squarefree strongly stable
by Lemma~\ref{lem:Tshift}.
For a squarefree strongly stable
ideal $J$, the graded Betti numbers
satisfy
$\beta_{i,i+j}(J)
=\sum_{\substack{v\in G(J)\\
\deg v=j}}
\binom{\max(v)-j}{i}
$
by \cite{AHH98}; see also
\cite[Section~7.4]{HH11}.
For $v=uw$, we have
$\max(v)=\max(u)$ and
$\deg v=\deg u+k-1$.
For each $u$, there are
$\binom{\min(u)-1}{k-1}$
choices of $w$.
This gives the Betti-number formula.
The regularity and projective-dimension
formulas follow immediately, using
$\reg I=\max\{\deg u:u\in G(I)\}$. 
Finally, if $I$ is the squarefree
Veronese ideal of degree $d+1$,
then $\sym{I}{k}$ is the squarefree
Veronese ideal of degree $d+k$,
provided the latter is nonzero.
\end{proof}

\begin{remark}\label{rem:beyondshifted}
The bound
$\reg\sym{I}{k}\le\reg I+k-1$
fails in general, even for ideals
with linear resolutions.
For example, let
$I=I(K_{2,2})=(x_1,x_2)(y_1,y_2)$.
Then $\reg I=2$, whereas
$\sym{I}{2}=(x_1x_2y_1y_2)$
has regularity $4$.
\end{remark}

\vspace*{1mm}
\noindent\textbf{Acknowledgments.}
 The author was partially supported by the
National Board for Higher Mathematics (NBHM) and the Science and
Engineering Research Board (SERB).

\vspace*{1mm}
\noindent\textbf{Data availability statement.}
Data sharing is not applicable to this article, as no datasets were
generated or analyzed during the current study.

\vspace*{1mm}
\noindent\textbf{Conflict of interest.}
The author declares that there are no competing interests or personal
relationships that could have influenced the work reported in this
paper.

\vspace*{1mm}
\noindent\textbf{Declaration of AI-assisted writing.}
The author used OpenAI's ChatGPT to assist in improving the language
and presentation of the manuscript. The author has reviewed
and edited the AI-assisted content and takes full responsibility for
the final version of the manuscript.

\bibliographystyle{abbrv}
\bibliography{ref.bib} 
\end{document}